\documentclass[12pt, a4paper, reqno, oneside]{article}

\usepackage{amssymb,amsmath}
\DeclareUnicodeCharacter{0335}{}
\usepackage{tikz-cd}
\usepackage{authblk} % For affiliations
\usepackage{physics} % For derivative notation (and more)
\usepackage{tikz}
\usepackage{titlesec} % Section/subsection/subsubsection style
\usetikzlibrary{babel}
\usepackage{adjustbox}
\usepackage{url}
\usepackage{xcolor}
\usepackage{yhmath}
\usepackage{eso-pic}
\usepackage[T1]{fontenc}
\usepackage[utf8]{inputenc}
\usepackage[english]{babel}
\usepackage{csquotes}
\usepackage{latexsym}
\usepackage{euscript}
\usepackage{mathtools}
\usepackage[all]{xy}
\usepackage{mathrsfs}
\usepackage[percent]{overpic}
\usepackage{setspace}
\usepackage{graphicx}
\usepackage{amsthm}
\usetikzlibrary{calc}
\usetikzlibrary{decorations.pathmorphing}

\tikzset{curve/.style={settings={#1},to path={(\tikztostart)
    .. controls ($(\tikztostart)!\pv{pos}!(\tikztotarget)!\pv{height}!270:(\tikztotarget)$)
    and ($(\tikztostart)!1-\pv{pos}!(\tikztotarget)!\pv{height}!270:(\tikztotarget)$)
    .. (\tikztotarget)\tikztonodes}},
    settings/.code={\tikzset{quiver/.cd,#1}
        \def\pv##1{\pgfkeysvalueof{/tikz/quiver/##1}}},
    quiver/.cd,pos/.initial=0.35,height/.initial=0}

\tikzset{tail reversed/.code={\pgfsetarrowsstart{tikzcd to}}}
\tikzset{2tail/.code={\pgfsetarrowsstart{Implies[reversed]}}}
\tikzset{2tail reversed/.code={\pgfsetarrowsstart{Implies}}}
\tikzset{no body/.style={/tikz/dash pattern=on 0 off 1mm}}
\usepackage[colorlinks=true,linktocpage=true, pagebackref=false, citecolor=black,linkcolor=black]{hyperref}
\usepackage{orcidlink}
\normalbaselines
\makeatletter

\def\ps@myfancy{\let\@mkboth\markboth
 \def\@evenhead{\vbox{\hsize\textwidth 
 \hbox to \textwidth{\sf\mdseries\thepage 
 \rule[-.6ex]{0mm}{2mm} \hfill\sf\large\leftmark}
 \vskip 1pt \hrule}}
 \def\@oddhead{\vbox{\hsize\textwidth 
 \hbox to \textwidth{{\sf\large\leftmark}
 \rule[-.6ex]{0mm}{2mm} \hfill\sf\mdseries{\thepage}}
 \vskip 1pt \hrule}}}

\def\ps@myfancyplain{
 \def\@evenhead{\vbox{\hsize\textwidth%
 \rule[-.6ex]{0mm}{2mm} \hfill }
 \vskip 1pt \hrule
 \vskip\headsep
 \vskip\textheight
 \vskip1pc
 \hbox to \textwidth{\sf\mdseries\thepage 
 \rule[-6ex]{0mm}{2mm} \hfill }}
 \def\@oddhead{\vbox{\hsize\textwidth 
 \vskip 1pt\hrule
 \vskip\headsep
 \vskip\textheight
 \vskip2pc
 \hbox to \textwidth{\hfill\rule[.4ex]{1pc}{2.5pt}
 \sf\mdseries\thepage}
}}}

\def\ps@myemptyfun{\ps@plain}

\makeatother
\providecommand{\proofname}{Demostraci\'on.}
 {\par\noindent{\it Demostraci\'on. }\nopagebreak\normalsize}%
 {}

 {\par\noindent{\it #1. }\nopagebreak\normalsize}%
 {\hfill\linebreak[2]\hspace*{\fill}$\square$\\[-1pt]}
\makeatother

\def\sqbullet{\raise.2ex\hbox{\vrule width 3.5pt height 3.5pt}}

{}

{\em}

\newcounter{substep}
\def\thesubstep{\arabic{substep}}

{\em}

\newcounter{subsubstep}
\def\thesubsubstep{\arabic{subsubstep}}

{\em}

\numberwithin{figure}{section}

{}

\newtheoremstyle{mystyle}%                % Name
  {}%                                     % Space above
  {}%                                     % Space below
  {\itshape}%                                     % Body font
  {}%                                     % Indent amount
  {\sf \bfseries}%                            % Theorem head font
  {}%                                    % Punctuation after theorem head
{ }%                                    % Space after theorem head, ' ', or \newline
  {\thmname{#1}\thmnumber{{\textcolor{blue}{\, \hspace{-1mm}#2.}}}\thmnote{ (#3)}}%                                     % Theorem head spec (can be left empty, meaning `normal')

\theoremstyle{mystyle}
\definecolor{royalblue(web)}{rgb}{0.25, 0.41, 0.88}
\hypersetup{
    colorlinks=true,
    linkcolor=royalblue(web),
    filecolor=magenta,      
    urlcolor=cyan,
    pdftitle={Overleaf Example},
    pdfpagemode=FullScreen,
    }

\titleformat{\section}
  {\normalfont\LARGE\bfseries\raggedright}{\thesection.}{1em}{}
\titleformat{\subsection}
  {\normalfont\Large\bfseries\raggedright}{\thesubsection.}{1em}{}
\titleformat{\subsubsection}
  {\normalfont\normalsize\itshape\raggedright}{\thesubsubsection.}{1em}{}

\newtheorem{Teor}{Theorem}[section]

\newtheorem{Prop}[Teor]{Proposition}
\newtheorem{Coro}[Teor]{Corollary}

\newtheorem{lema}[Teor]{Lema}
\newtheorem{Defi}[Teor]{Definition}

\newtheorem{Obse}[Teor]{Remark}

 \newcommand{\N}{{\mathbb N}}
\newcommand{\Z}{{\mathbb Z}} 
\newcommand{\R}{{\mathbb R}}

\newcommand{\m}{\operatorname{mult}}

\newcommand{\mail}[1]{\small\href{mailto:#1}{#1}}

\newenvironment{Abstract}
{
\begin{center}
\textbf{Abstract}\\
\vspace{0.25cm}
\begin{minipage}{14.5cm}}
{\footnotesize
\end{minipage}
\end{center}}
\begin{document}
%-------------------------------------------
% Title, authors and affiliations
%-------------------------------------------

%\maketitle

	\begin{center}
		{\huge {\bfseries 
        Interpolation of nonlocal Sobolev spaces for general gradients}\par}
		\vspace{1cm}
		
\begin{center} 
	 {\Large Guillermo García-Sáez}{\small\textsuperscript{1}} \\
	\mail{guillermo.garciasaez@uclm.es}\:\orcidlink{0009-0008-4335-119X}
\end{center}
\vspace{5mm}

\textsc{\textsuperscript{1}ETSII, Departamento de Matem\'aticas\\ Universidad de Castilla-La Mancha} \\
		Campus Universitario s/n, 13071 Ciudad Real, Spain. \\ \vspace{5mm}
\end{center}
\begin{Abstract}
We study the interpolation structure of nonlocal Sobolev spaces associated with general radial gradients. For compactly supported kernels satisfying suitable regularity and two-sided growth assumptions of fractional order, we characterize the spaces $H^{\rho,p}(\mathbb{R}^n)$, $1<p<\infty$, by the Gustavsson--Peetre interpolation method between $L^p(\mathbb{R}^n)$ and $W^{1,p}(\mathbb{R}^n)$. The function parameter is determined, up to equivalence, by the behavior of the kernel near the origin. We also establish a complex interpolation formula between spaces associated with different kernels and integrability exponents, in which the interpolated kernel is their geometric mean. For the corresponding nonlocal Gagliardo spaces $W^{\rho,p}(\R^n)$, defined through an integral seminorm of first differences, we prove an identification with generalized Besov spaces and a real interpolation characterization involving the same function parameter. These results extend the classical relation between Bessel potential and Gagliardo spaces beyond power-type kernels, yielding continuous embeddings between the two families and their coincidence in the Hilbertian case. We further discuss logarithmic kernels and restriction spaces on Lipschitz domains, and obtain extension and compact embedding results for the Gagliardo spaces under weaker assumptions on the kernel.
\end{Abstract}
\tableofcontents
\section{Introduction}

Fractional Sobolev spaces provide a natural framework for measuring
regularity between Lebesgue spaces and first-order Sobolev spaces.
Two of their main realizations are the Bessel potential spaces
$H^{s,p}(\mathbb{R}^n)$ and the Gagliardo spaces
$W^{s,p}(\mathbb{R}^n)$, where $0<s<1$ and $1<p<\infty$.
Both are interpolation spaces between $L^p(\mathbb{R}^n)$ and
$W^{1,p}(\mathbb{R}^n)$: the former arise from complex interpolation,
whereas the latter are obtained by the real method with exponent $p$;
see \cite{BellidoGarciaI2025,BerghLofstrom1976,Triebel1995}.
These descriptions explain their common regularity scale, their
different behavior away from the Hilbertian setting, and many of
their embedding properties. For nonlocal spaces governed by a
general interaction kernel, however, a single differentiability
index need not describe the regularity. This raises a natural
question: how does the kernel determine the interpolation structure,
and what becomes of the relation between the Bessel and Gagliardo
constructions?

The connection of Bessel potential spaces with nonlocal differential operators is provided by
the Riesz fractional gradient, a fractional integro-differential operator defined up to a normalization constant as $$D^s u(x)=\int_{\R^n}\frac{u(x)-u(y)}{|x-y|^{n+s}}\frac{x-y}{|x-y|}\,dy,\,x\in \R^n,\,u\in C_c^\infty(\R^n),$$ introduced by Shieh and Spector in their celebrated works \cite{ShiehSpector2015,ShiehSpector2018}. Among many other relevant results, the authors showed that the
associated Sobolev spaces coincide with the Bessel potential spaces
with equivalence of norms.
This allows fractional regularity to be described by the
integrability of a vector-valued nonlocal derivative, and provides
a framework for variational problems involving fractional gradients.
Subsequent developments include gradients with finite interaction
range and more general radial kernels
\cite{BellidoCuetoMora2023, BMS, Cueto2023, CuetoHidde2025, Kreisbeck2024}. Such operators have found many relevant applications like, among others, in
nonlocal models of elasticity \cite{BellidoCuetoMoraII2023, Bellido_NH2024, BellidoGarcia2026, BellidoGarcia2026II, Hidde2026, KreisbeckMoraSeifert}.
In the general setting considered here, which is the one developed in \cite{BMS}, the nonlocal gradient associated
with a radial kernel $\rho$ is defined on smooth functions of
compact support by
\[
    D_\rho u(x)
    =\int_{\mathbb{R}^n}
      \frac{u(x)-u(y)}{|x-y|}\,
      \frac{x-y}{|x-y|}\,\rho(x-y)\,dy,
\]
and the corresponding space $H^{\rho,p}(\mathbb{R}^n)$ consists of
the functions $u\in L^p(\mathbb{R}^n)$ whose weak nonlocal gradient
belongs to $L^p(\mathbb{R}^n;\mathbb{R}^n)$.
The choice $\rho(x)=c_{n,s}|x|^{-n-s+1}$ recovers the Riesz
fractional gradient and hence $H^{s,p}(\mathbb{R}^n)$.

The purpose of this work is to study these spaces from the point
of view of interpolation theory and to relate them to the
corresponding Gagliardo spaces based on general kernels recently introduced in \cite{dorer}. We mainly consider compactly
supported radial kernels satisfying the assumptions introduced
in \cite{BMS}, with lower and upper growth bounds of fractional
order. These bounds allow the singularity of the kernel to depart
from a pure power, as happens, for example, for kernels with
logarithmic factors. The appropriate interpolation parameter must
therefore retain more information than a constant exponent.
To this end, we use the Gustavsson--Peetre method
\cite{GustavssonPeetre} and real interpolation with a function
parameter \cite{Gustavsson}, together with the theory of spaces
of generalized smoothness \cite{AnopMurach, Farkas}.

A starting point is the recent Fourier characterization of
nonlocal Sobolev spaces obtained by Happ, Kreisbeck and Mora-Corral
\cite{Happ2026}. Their result identifies
$H^{\rho,p}(\mathbb{R}^n)$ through a Fourier weight associated with
$D_\rho$. Building on this characterization, we construct a smooth
weight $\alpha_\rho$ of generalized smoothness and an admissible
interpolation parameter $\varphi_\rho$. If $\overline\rho$ denotes
the radial representation of $\rho$, the weight satisfies
$\alpha_\rho(t)\asymp t^{1-n}\overline\rho(1/t)$ for sufficiently
large $t$. We set $\varphi_\rho(t)=\alpha_\rho(t)$ for $t\geq1$,
and $\varphi_\rho(t)=1/\alpha_\rho(1/t)$ for $0<t<1$, and prove, with equivalence of norms, that
\[
    H^{\rho,p}(\mathbb{R}^n)
    =H^{\alpha_\rho(\cdot),p}(\mathbb{R}^n)
    =\big\langle L^p(\mathbb{R}^n),
                 W^{1,p}(\mathbb{R}^n)\big\rangle_{\varphi_\rho},
    \qquad 1<p<\infty,
\]
where the brackets denote the Gustavsson--Peetre interpolation
method. Thus, the behavior of the interaction kernel near the
origin determines the interpolation parameter up to equivalence
and places the nonlocal space within the scale of generalized
Bessel potential spaces. For a fractional kernel, the parameter
has the familiar power behavior $t^s$ near zero.

This identification also allows us to interpolate between spaces
associated with different kernels and different integrability
exponents. More precisely, for admissible kernels $\rho_0,\rho_1$,
$1<p_0,p_1<\infty$ and $0<\theta<1$, we obtain
\[
    \begin{gathered}
        \big[H^{\rho_0,p_0}(\mathbb{R}^n),
             H^{\rho_1,p_1}(\mathbb{R}^n)\big]_\theta
        =H^{\rho_\theta,p_\theta}(\mathbb{R}^n),
        \\
        \rho_\theta=\rho_0^{1-\theta}\rho_1^\theta,
        \qquad
        \frac1{p_\theta}
        =\frac{1-\theta}{p_0}+\frac{\theta}{p_1},
    \end{gathered}
\]
again with equivalence of norms. In particular, complex
interpolation has a direct description in terms of the interaction
kernels: it corresponds to their geometric mean. This extends the
interpolation law for fractional differentiability orders to kernels
whose singularities are not described by a single power.

The second family considered in this work is defined directly
through differences. We consider the nonlocal Gagliardo spaces
studied in \cite{carioni, dorer}: the spaces
$W^{\rho,p}(\mathbb{R}^n)$ defined as the spaces of functions $u\in L^p(\mathbb{R}^n)$ such that the nonlocal seminorm
\[
    [u]_{\rho,p}^{\,p}
    :=\int_{\mathbb{R}^n}\int_{\mathbb{R}^n}
       \frac{|u(x)-u(y)|^p}{|x-y|^{\,n-p(n-1)}}
       \rho(x-y)^p\,dx\,dy,
\] is finite.
For the fractional kernel, this gives the usual Gagliardo seminorm,
up to its normalization. Under the same assumptions on $\rho$,
we prove that
\[
    W^{\rho,p}(\mathbb{R}^n)
    =B^{\alpha_\rho(\cdot)}_{p,p}(\mathbb{R}^n)
    =\big(L^p(\mathbb{R}^n),
          W^{1,p}(\mathbb{R}^n)\big)_{\varphi_\rho,p},
    \qquad 1<p<\infty,
\] which places the nonlocal Gagliardo spaces within the scale of generalized smoothness Besov spaces via the generalized real interpolation method.
Consequently, the same function parameter describes both nonlocal
constructions, while the choice of interpolation method
distinguishes them. This yields continuous embeddings from
$W^{\rho,p}$ into $H^{\rho,p}$ for $1<p\leq2$, and in the reverse
direction for $2\leq p<\infty$. In particular, the two spaces
coincide when $p=2$, with equivalence of norms, as in the classical
fractional setting. In particular, we will also establish the equivalence between both $H^{\rho,2}$ and $W^{\rho,2}$ without interpolation theory proving the equivalence $$1+m_{\rho}(\xi)\asymp 1+|\lambda_\rho(\xi)|^2,$$ where $m_\rho$ is the symbol that appears in the Fourier characterization of the Gagliardo seminorm for $p=2$, and $\lambda_\rho(\xi)$ is the Fourier symbol of the nonlocal gradient $D_\rho$.

We also consider restriction spaces on Lipschitz domains and
illustrate the potential description with logarithmic kernels,
relating the resulting spaces to logarithmic Bessel potential
spaces \cite{Chen, OpicTrebels}. Finally, we study some structural properties
of the Gagliardo spaces that remain available under weaker assumptions
than those used for the interpolation characterizations.
By verifying the hypotheses of the integral-kernel theory in
\cite{Brusca}, we obtain extension and compact embedding results
on bounded Lipschitz domains, together with
Poincar\'e--Wirtinger inequalities when the domain is additionally
connected.

The paper is organized as follows. Section~2 fixes the notation.
Section~3 recalls the nonlocal gradient setting and its potential
characterization, including the logarithmic example. Section~4
collects the interpolation results and the spaces of generalized
smoothness used throughout the paper. Section~5 establishes the
interpolation description of $H^{\rho,p}$, the complex interpolation
formula and the corresponding results for restriction spaces on
domains. Finally, Section~6 introduces the nonlocal Gagliardo
spaces, proves their real interpolation characterization and
studies their structural properties, the Hilbertian case and
how they relate to a broader class of kernels.
\section{Notation}
We fix $n\in \mathbb{N}$ the dimension of our ambient space $\R^n$ and we will denote by $\Omega\subset \R^n$ an open bounded subset representing the body. The notation for Sobolev $W^{1,p}$ and Lebesgue $L^p$ spaces is the standard one, as is that of smooth
	functions of compact support $C_c^\infty$. We will
	indicate the domain of the functions, as in $L^p
	(\Omega)$; the target is indicated only if it is not $\R$. The notation for the classical gradient in its weak sense is the usual one $D$. For $\alpha \in\mathbb{N}^n$, we give the standard meaning to the partial derivative $\partial^\alpha$ and the size $|\alpha|$. For derivatives of a function of one variable, we use the notation $f^{k)}(t)=\frac{d^kf}{dt^k}$. For a Banach space $X$, we denote its topological dual as $X^*$, and the usual {\it strong} convergence in normed spaces is denoted as $\to$.

	We write $|x|=\left(\sum_{j=1}^n x_j^2\right)^{1/2}$ for the Euclidean norm of a vector $x=(x_1,\ldots,x_n)\in \R^n$. The ball centered at $x\in \R^n$ and with radius $r>0$ is denoted by $B_r(x)=\{y\in \R^n: |x-y|<r\}$ or if there is no ambiguity about the center, just by $B_r$. The complement of a set $E\subset \R^n$ is denoted by $E^c=\R^n\setminus{E}$, its closure by $\overline{E}$ and its boundary by $\partial E$. We will denote that a subset $K$ of a set $E$ is a compact subset relative to the topology of $E$ by $K\subset \subset E$. We denote by $\omega_{n-1}$ the surface measure of the $(n-1)$-dimensional sphere $\mathbb{S}^{n-1}$.
	
	Our convention for the Fourier transform of functions $f\in L^1(\R^n)$ is $$\widehat{f}(\xi)=\int_{\R^n}f(x)e^{-2\pi ix\cdot \xi}\,dx,\,\xi\in \R^n.$$ This definition is extended by continuity and duality to other function and distribution spaces as usual in the theory of function spaces. We will sometimes use the alternative notation $\mathfrak{F}(f)$ for $\widehat{f}$. More details about this operator can be found in the classical text \cite{Grafakos2008}. We denote by $\mathcal{S}$ the Schwartz space and $\mathcal{S}'$ the space of tempered distributions.
	
	Regarding radial functions we have the following definitions:
	\begin{itemize}
		\item A function $f:\R^n\to \R$ is \textit{radial} if there exists $\overline{f}:[0,\infty)\to \R$ such that $\overline{f}(|x|)=f(x)$ for every $x\in \R^n$. The function $\overline{f}$ is called the \textit{radial representation} of $f$.
		\item A radial function $f:\R^n\to \R$ is \textit{radially decreasing} if its radial representation is a
		decreasing function.
		\item A function $g:\R^n\to \R^n$ is \textit{vector radial} if there exists a radial function $\overline{\phi}:[0,\infty)\to \R$  such that
		$g(x)=\overline{\phi}(|x|)x$ for every $x\in \R^n$.
	\end{itemize}
	
	For real valued functions, we use the monotonicity properties of being increasing and decreasing in
	the non-strict sense. A function $f :\R \to\R$ is called \textit{almost decreasing} if there exists a positive constant $C$ such
	that $f(t)\geq Cf(s)$ for every $t\leq s$, and an analogous definition for \textit{almost increasing}.

    Let $(E_0,E_1)$ be a couple of Banach spaces. We say that the couple is \textit{compatible} if there exists a Hausdorff topological vector space $\mathcal{E}$ such that $E_0,E_1\xhookrightarrow{}\mathcal{E}$. We say that a Banach space $E$ is intermediate with respect to the couple if $E_0\cap E_1\xhookrightarrow{}E\xhookrightarrow{}E_0+E_1$. Let $(F_0,F_1)$ be another compatible couple and let $T:E_0+E_1\to F_0+F_1$ be a bounded linear operator. We say that $T$ is \textit{admissible} if $T:E_j\to F_j$, $j=0,1,$ continuously. Given $E$ an intermediate space with respect to the couple $(E_0,E_1)$, and $F$ an intermediate space with respect to $(F_0,F_1)$, we say that $E,F$ are \textit{interpolation spaces} with respect to the couples $(E_0,E_1)$ and $(F_0,F_1)$, respectively, if for every admissible operator $T:E_0+E_1\to F_0+F_1$, $T:E\to F$ continuously. The methods to construct such interpolation spaces for given couples are called \textit{interpolation functors.} We refer to \cite{BerghLofstrom1976,GarciaSaez2024,Lunardi2018,Triebel1995} for detailed expositions on the interpolation theory and the most important interpolation functors.
\section{Nonlocal Sobolev spaces and potential theory}

The kernels considered here are the extension of the fractional-gradient setting of
Shieh and Spector \cite{ShiehSpector2015, ShiehSpector2018} done by 
Bellido, Mora-Corral and Sch\"onberger in \cite{BMS} (see also \cite{BellidoGarciaNiño2026, CuetoHidde2025}).

Fix $n\geq1$ and $1<p<\infty$. 
For a radial kernel $\rho(x)=\overline{\rho}(|x|)$, the nonlocal gradient $D_\rho$ is
defined on test functions by
\[
 D_\rho u(x)=\int_{\R^n}
 \frac{u(x)-u(y)}{|x-y|}\frac{x-y}{|x-y|}\rho(x-y)\,dy,\,x\in \R^n.
\]

For the kernels $\rho$, we will use the hypotheses of \cite{BMS}. We write $f_\rho(t):=t^{n-2}\overline{\rho}(t)$. 
\begin{enumerate}
\item[(H0)] $\rho:\R^n\setminus\{0\}\to[0,\infty)$ is radial and finite,
$\rho\in L^1_{\mathrm{loc}}(\R^n)$, and
\[
 \int_{\R^n}\min\{1,|x|^{-1}\}\rho(x)\,dx<\infty,
 \qquad \inf_{x\in B_\varepsilon(0)\setminus\{0\}}\rho(x)>0,
\] for some $\varepsilon>0$.
\item[(H1)] $f_\rho$ is nonincreasing on $(0,\infty)$, and there exists
$0<\mu<1$ such that
$f_\rho(t)\leq\mu f_\rho(t/2)$ for $0<t<\varepsilon$.
\item[(H2)] $f_\rho\in C^\infty(0,\infty)$, and for $0<t<\varepsilon$,
\[
 -f_\rho'(t)\geq c\,t^{-1}f_\rho(t),\qquad
 |f_\rho^{(k)}(t)|\leq C_k t^{-k}f_\rho(t)\quad(k\geq1),
 \qquad c>0.
\]
\item[(H3)] For some $\sigma\in(0,1)$,
$t\mapsto t^{n+\sigma-1}\bar\rho(t)$ is almost nonincreasing on
$(0,\varepsilon)$.
\item[(H4)] For some $\gamma\in(0,1)$,
$t\mapsto t^{n+\gamma-1}\bar\rho(t)$ is almost nondecreasing on
$(0,\varepsilon)$.
\end{enumerate}

In \cite{BMS, CuetoHidde2025}, examples of admissible kernels are given, the most paradigmatic one 
	being $$\rho^s(x):=c_{n,s}\frac{1}{|x|^{n+s-1}},\,s\in (0,1),$$ where $c_{n,s}$ is a suitable normalization constant. The nonlocal gradient associated with this kernel is the well-known Riesz fractional gradient $D^s=:D_{\rho^s}$, which is defined for functions $u\in C_c^\infty(\R^n)$ as $$D^su(x)=c_{n,s}\int_{\R^n}\frac{u(x)-u(y)}{|x-y|^{n+s}}\frac{x-y}{|x-y|}\,dy=D(I_{1-s}u),$$ where $I_s$ is the Riesz potential.
	The Riesz fractional gradient was introduced by Shieh and Spector in \cite{ShiehSpector2015,ShiehSpector2018}, and since then has been widely studied.
	
	Another important example concerns a truncated version of the Riesz fractional gradient introduced in \cite{BellidoCuetoMora2023}, which is suitable for bounded domains. Let $\delta>0$ and let $w_\delta$ be a non-negative radial cut-off
	function such that $w_\delta(0)>0$, $w_\delta\in C_c^\infty\left(B(0,\delta)\right)$ and $w_\delta/|\cdot|^{1+s}$ is radially decreasing.
	Now, taking the kernel 
	$$\rho_\delta^s(x):=\frac{w_\delta(x)}{|x|^{n+s-1}},$$ 
	yields that $D_{\rho_\delta^s}$ is, up to a normalization constant, equal to $D_\delta^s$, the nonlocal gradient introduced in \cite{BellidoCuetoMora2023}.

More interesting examples are the Riesz fractional gradient with (radial) variable differentiability parameter given by $$\rho^{s(\cdot)}(x):=\frac{\chi(x)}{|x|^{n+s(|x|)-1}},$$ where $\chi$ is a nonnegative radial function in $C_c^\infty(\R^n)$ such that $\chi(0)>0$ and $\chi(x)/|x|^{1+s(|x|)}$ is radially decreasing, and $s:[0,\infty)\to (0,1)$ is smooth. In this case, one can take $\sigma=\min_{[0,\varepsilon]}s$ and $\gamma=\max_{[0,\varepsilon]}s$, for any $\varepsilon>0$. Another example is the logarithmic kernel
\[
   \rho_{\log}^s(x)
   =\frac{-\chi(x)\log(|x|)}{|x|^{n+s-1}},\qquad x\ne0,
\]
where $0<s<1$ and $\chi\in C_c^\infty(B_1)$ is nonnegative, radial,
nonincreasing in the radial variable and equal to one near the origin.
In this case one may take $\sigma=s$ and any $\gamma\in(s,1)$.

One of the most important features of nonlocal gradients is that they can be written as a convolution of a classical gradient with a potential. Namely, if we define $$Q_\rho(x)=\int_{|x|}^\infty\frac{\bar\rho(t)}t\,dt,\,x\in\R^n\setminus\{0\},$$ we have that $D_\rho u=Q_\rho*D u=D(Q_\rho*u)$, for functions $u\in C_c^\infty(\R^n)$ from
\cite[Proposition~2.6]{BMS}, and hence 
$\widehat{D_\rho u}(\xi)=2\pi i\xi q_\rho(\xi)\widehat u(\xi)$, where $q_\rho=\widehat{Q_\rho}$. Moreover, if $\rho\in L^1(\R^n)$, then $Q_\rho\in L^1(\R^n)$ (in fact, $\rho$ having compact support is enough for $Q_\rho$ to lie in $L^1(\R^n)$, \cite[Lemma 2.5 (iii)]{BMS}). For our development, it will be crucial to control the asymptotics of $q_\rho$. In particular we will require the following result from \cite[Lemma 4.3, Lemma 4.10, Lemma 7.1]{BMS}.
\begin{lema}\label{estimacionesQ}Suppose that $\rho$ satisfies $(H0)-(H4)$ with compact support. Then, $\widehat{Q}_\rho$ is positive and there exists a positive constant $c$ such that $$q_\rho(\xi)\geq c|\xi|^{-n}\overline{\rho}(1/|\xi|),\,\xi \in B_{2/\varepsilon}^c.$$ Furthermore, if $\rho$ satisfies $(H4)$, then there exists a positive constant $C$ such that $$q_\rho(\xi)\leq C|\xi|^{-n}\overline{\rho}(1/|\xi|),\,\xi\in B_{1/\varepsilon}^c.$$ Moreover, for every multi-index $\beta$, $$|\partial^\beta q_\rho(\xi)|\leq C_\beta \langle\xi\rangle^{-|\beta|}q_\rho(\xi).$$
    
\end{lema}

For a vector field $v\in C_c^\infty(\R^n;\R^n)$, we can define its \textit{nonlocal divergence} as $$\operatorname{div}_\rho v(x):=\int_{\R^n}\frac{v(x)-v(y)}{|x-y|}\cdot \frac{x-y}{|x-y|}\rho(x-y)\,dy,\,x\in \R^n.$$ The nonlocal divergence is the dual operator of the nonlocal
gradient in the sense of integration by parts by \cite[Proposition 3.2]{BMS}, i.e., for $u\in C_c^\infty(\R^n)$, we have that $$\int_{\R^n} D_\rho u(x)\cdot v(x)\,dx=-\int_{\R^n}u(x)\operatorname{div}_\rho v(x)\,dx,$$ for every $v\in C_c^\infty(\R^n;\R^n)$. Moreover, we can define the nonlocal gradient distributionally for $u\in L^1(\R^n)+L^\infty(\R^n)$ as the distribution $\mathfrak{G}_\rho u$ such that $$\langle\mathfrak{G}_\rho u,v\rangle=-\int_{\R^n}u(x)\operatorname{div}_\rho v(x)\,dx,\,v\in C_c^\infty(\R^n;\R^n),$$ which coincides with the previous definition when $$\int_K\int_{\R^n}\frac{|u(x)-u(y)|}{|x-y|}\rho(x-y)\,dx\,dy<\infty,\,\forall K\subset\subset \R^n.$$

Using its distributional extension, \textit{nonlocal Sobolev spaces} for a general kernel $\rho$ were introduced in \cite{BMS} as the space
\[
 H^{\rho,p}(\R^n)=\{u\in L^p(\R^n):\mathfrak{G}_\rho u\in L^p(\R^n;\R^n)\},
 \qquad \|u\|_{H^{\rho,p}}=\|u\|_p+\|\mathfrak{G}_\rho u\|_p.
\] For an open set $\Omega\subset \R^n$, we can also define the closed subspace $$H^{\rho,p}_0(\Omega):=\{u\in H^{\rho,p}(\R^n): u|_{\Omega^c}=0\}.$$
We recall that the choice $\rho=\rho^s$ leads to the Bessel potential spaces $H^{s,p}$ (see \cite{ShiehSpector2015}). By \cite[Theorem 3.9]{BMS}, under hypothesis $(H0)$, we have the density results $H^{\rho,p}(\R^n)=\overline{C_c^\infty(\R^n)}^{\norm{\cdot}_{\rho,p}}$, and $H^{\rho,p}_0(\Omega)=\overline{C_c^\infty(\Omega)}^{\norm{\cdot}_{\rho,p}}$, with the last being true if $\Omega$ is Lipschitz.

An interesting observation is that we can restrict ourselves to kernels with compact support without loss of generality. Indeed, choose $\chi\in C^\infty([0,\infty))$ with $0\leq\chi\leq1$,
nonincreasing, equal to one near zero, and compactly supported.
The kernel $\widetilde\rho(x)=\chi(|x|)\rho(x)$ still satisfies
(H0)--(H4).
Moreover,
$(\rho-\widetilde\rho)/|\cdot|\in L^1(\R^n)$ by (H0), and by Young's
inequality,
\begin{equation}\label{eq:tail}
 \|D_\rho u-D_{\widetilde\rho}u\|_p
 \leq2\left\|\frac{\rho-\widetilde\rho}{|\cdot|}\right\|_1\|u\|_p.
\end{equation}
It follows that
$H^{\rho,p}=H^{\widetilde\rho,p}$ with equivalence of norms by
\cite[Proposition~3.10, p.~15]{BMS}.
\subsection{Potential characterization}
In order to study the nonlocal spaces $H^{\rho,p}(\R^n)$ as interpolation spaces, as it happened for the purely fractional case with $\rho=\rho^s$, we would require a characterization of the spaces in terms of a potential. We recall that the \textit{Bessel potential spaces} $H^{s,p}(\R^n)$ (see \cite{BellidoGarciaI2025} and the references therein) are defined for $s \in \R$ and $p \in (1,\infty)$ by
\[
H^{s,p}(\R^n):=\left\{u \in \mathcal{S}'(\R^n)\,:\, \left(\langle \cdot \rangle^{s}\widehat{u}\right)^{\vee} \in L^p(\R^n)\right\} \quad \text{with} \quad \norm{u}_{H^{s,p}(\R^n)}:=\norm{\left(\langle \cdot \rangle^{s}\widehat{u}\right)^{\vee}}_{L^p(\R^n)},
\]
where $\langle \xi \rangle :=(1+|\xi|^2)^{1/2}$, sometimes referred to as the \textit{Japanese bracket}, is a smooth function that behaves like $|\xi|$ for large $\xi$. For $0<s<1$, these spaces are a particular example of the nonlocal Sobolev spaces $H^{\rho,p}(\R^n)$ for the choice $\rho=\rho^s$ as it was proved by Shieh and Spector in \cite{ShiehSpector2015}. The operator $$\Lambda_s:\mathcal{S}'(\R^n)\to \mathcal{S}'(\R^n); f\mapsto \left(\langle \cdot \rangle^{-s}\widehat{f}\right)^{\vee},$$ is the \textit{Bessel potential} of order $s$, and for $s>0$ it can be shown that there exists a function $g_s\in L^1(\R^n)$ such that $\Lambda_sf=g_s*f$. This identification and the properties of the potential are the key to establishing that the spaces $H^{s,p}(\R^n)$ can be obtained by complex interpolation \cite[Theorem 3.7]{BellidoGarciaI2025}. In particular, for every $s>0$, $$H^{s,p}(\R^n)=[L^p(\R^n),W^{k,p}(\R^n)]_{\theta},$$ for $\theta\in (0,1)$ and $k\in \mathbb{N}$ such that $k\theta=s$. In order to obtain such a result for some interpolation method, we would need some characterization of the spaces $H^{\rho,p}(\R^n)$ analogous to the fractional one. This is the main objective of the recent work \cite{Happ2026}, where several characterizations for the spaces $H^{\rho,p}(\R^n)$, analogous to the classical ones, are established. 

The easiest case is the Hilbertian one, where no multiplier theory is needed. In \cite{BellidoGarciaNiño2026}, it was established that for the case $p=2$, we could characterize the space $H^{\rho,2}(\R^n)$ as the potential space $$H^{\rho,2}(\R^n)=\{u\in \mathcal{S}'(\R^n):w_\rho(\xi)\widehat{u}(\xi)\in L^2(\R^n)\},$$ where
$$
w_\rho(\xi):=\bigl(1+4\pi^2|\xi|^2q_\rho(\xi)^2\bigr)^{1/2}=\left(1+|\lambda_\rho(\xi)|^2\right)^{1/2},\,\xi\in \R^n.$$ This result was extended to $1<p<\infty$, by Happ, Kreisbeck and Mora-Corral in \cite[Theorem 3.1]{Happ2026}, having that
\begin{equation*}
 H^{\rho,p}(\R^n)=\{u\in\mathcal{S}':\mathcal F^{-1}(w_\rho \widehat u)\in L^p\},
 \qquad \|u\|_{H^{\rho,p}}\asymp\norm{\mathfrak{F}^{-1}(w_\rho \widehat{u})}_p,\qquad 1<p<\infty,
\end{equation*}
with equivalence of the norms. This identification will be the key for our development as well as the following technical lemma:

\begin{lema}\label{lemma-clave}Suppose that $\rho$ is a compactly supported kernel satisfying $(H0)-(H4)$. Then, there exists $R\geq 1$ such that $$w_\rho(\xi)\asymp |\xi|^{1-n}\overline{\rho}(1/|\xi|),\,|\xi|\geq R.$$ Moreover, for every multi-index $\beta$, $$|\partial^\beta w_\rho(\xi)|\leq C_\beta \langle\xi\rangle^{-|\beta|}w_\rho(\xi),\,\xi\in \R^n.$$
\end{lema}
\noindent{}\textbf{Proof:}
First of all, note that for a large enough $|\xi|$, $w_\rho(\xi)=\left(1+4\pi^2|\xi|^2q_\rho(\xi)^2\right)^{1/2}=\langle 2\pi\xi q_\rho(\xi)\rangle\asymp |\xi|q_\rho(\xi)$. Now, for $|\xi|\geq 2/\varepsilon$, we have that by Lemma \ref{estimacionesQ}, $$q_\rho(\xi)\asymp |\xi|^{-n}\overline{\rho}(1/|\xi|)\implies |\xi|q_\rho(\xi)\asymp |\xi|^{1-n}\overline{\rho}(1/|\xi|),$$ and hence there exists $R\geq 2/\varepsilon$ such that $$w_\rho(\xi)\asymp |\xi|q_\rho(\xi)\asymp |\xi|^{1-n}\overline{\rho}(1/|\xi|),\,|\xi|\geq R.$$ The derivative estimates follow from
\cite[Lemma~2.8, estimate~(2.24)]{Happ2026}, applied with exponent
$1/2$. Indeed,
we obtain, for every multi-index $\beta$ and $|\xi|\geq 1$,
\[
|\partial^\beta w_\rho(\xi)|
\leq C_\beta |\xi|^{-|\beta|}|\lambda_\rho(\xi)|
\leq C_\beta |\xi|^{-|\beta|}w_\rho(\xi)
\leq \widetilde C_\beta
\langle\xi\rangle^{-|\beta|}w_\rho(\xi).
\]
On the other hand, $Q_\rho\in L^1(\mathbb R^n)$ has compact
support, so $q_\rho$ and $w_\rho$ are smooth on $\mathbb R^n$.
Thus, setting
\[
M_\beta:=\sup_{|\xi|\leq1}|\partial^\beta w_\rho(\xi)|<\infty
\]
and using $w_\rho\geq1$, we have
\[
|\partial^\beta w_\rho(\xi)|
\leq M_\beta
\leq 2^{|\beta|/2}M_\beta
\langle\xi\rangle^{-|\beta|}w_\rho(\xi),
\qquad |\xi|\leq1.
\]
Combining both estimates proves the claimed bound on
$\mathbb R^n$. \qed

\subsection{The scale of logarithmic Bessel potential spaces}
Recently, the scale of logarithmic Bessel potential spaces $H^{s,p,\lambda}_{\log}(\R^n)$, defined for $s>0$, $\lambda>1$ and $1<p<\infty$ as $$H^{s,p,\lambda}_{\log}(\R^n):=\{u\in \mathcal{S}'(\R^n):\mathfrak{F}^{-1}\left((\lambda+4\pi^2|\xi|^2)^{s/2}\log (\lambda+4\pi^2|\xi|^2)\widehat{u}(\xi)\right)\in L^p(\R^n)\},$$
has been extensively studied in \cite{Chen}. The definition of those spaces does not depend on the value $\lambda$ \cite[Proposition 1.10]{Chen}, but they are strictly embedded into the critical case $\lambda=1$, \cite[Proposition 1.11]{Chen}. A characterization without the parameter $\lambda$ is given in \cite[Proposition 1.13]{Chen}, where they are identified with the scale of Opic and Trebels \cite{OpicTrebels}, given by $$H^{s,p}_{\log}(\R^n):=\{u\in \mathcal{S}'(\R^n): \mathfrak{F}^{-1}\left((1+4\pi^2|\xi|^2)^{s/2}(1+\log (1+4\pi^2|\xi|^2))\widehat{u}(\xi)\right)\in L^p(\R^n)\}$$ Moreover, they are intermediate with respect to the Bessel scale \cite[Proposition 1.12]{Chen}, i.e., $$H^{s+\delta,p}(\R^n)\xhookrightarrow{}H^{s,p,\lambda}_{\log}(\R^n)\xhookrightarrow{}H^{s,p}(\R^n),\,\delta>0.$$
It can be shown that in fact the logarithmic Bessel scale can be written in terms of nonlocal Sobolev spaces for a particular choice of kernel $\rho$. In particular, for the choice $\rho=\rho^s_{\log}$, we have the following result:
\begin{Teor}
    Let $s\in (0,1)$, $1<p<\infty$ and $\lambda>1$. Then, we have that $$H^{\rho^s_{\log},p}(\R^n)=H^{s,p,\lambda}_{\log}(\R^n),$$ with equivalence of the norms.
\end{Teor}
\noindent{}\textbf{Proof:} 
We take the cut-off function $\chi$ defining $\rho_{\log}^s$
to be radial, nonnegative and nonincreasing, with
$\chi\in C_c^\infty(B_1)$ and $\chi=1$ in a neighborhood
of the origin. Then $\rho_{\log}^s$ satisfies (H0)--(H4),
with $\sigma=s$ and any $\gamma\in(s,1)$, after choosing
the radius in these hypotheses sufficiently small.
Set
\[
w(\xi):=w_{\rho_{\log}^s}(\xi),
\qquad
m(\xi):=(\lambda+4\pi^2|\xi|^2)^{s/2}
\log(\lambda+4\pi^2|\xi|^2).
\]
By Lemma \ref{lemma-clave} and the choice of $\chi$, for $|\xi|$
sufficiently large,
\[
w(\xi)
\asymp |\xi|^{1-n}\bar\rho_{\log}^s(1/|\xi|)
=|\xi|^s\log|\xi|
\asymp m(\xi).
\]
Moreover, both symbols are smooth and strictly positive
on $\mathbb R^n$, since $w\geq1$ and
$m(\xi)\geq\lambda^{s/2}\log\lambda>0$.
Their comparison therefore extends to all frequencies:
\[
c\,m(\xi)\leq w(\xi)\leq C\,m(\xi),
\qquad \xi\in\mathbb R^n,
\]
for some constants $c,C>0$.
Lemma \ref{lemma-clave} and direct differentiation of $m$ give,
for every multi-index $\beta$,
\[
|\partial^\beta w(\xi)|
\leq C_\beta\langle\xi\rangle^{-|\beta|}w(\xi),
\qquad
|\partial^\beta m(\xi)|
\leq C_\beta\langle\xi\rangle^{-|\beta|}m(\xi).
\]
Differentiating the identities $w\,w^{-1}=1$ and
$m\,m^{-1}=1$, induction on $|\beta|$ yields the
corresponding estimates for the reciprocals:
\[
|\partial^\beta(w^{-1})(\xi)|
\leq C_\beta\langle\xi\rangle^{-|\beta|}w(\xi)^{-1},
\qquad
|\partial^\beta(m^{-1})(\xi)|
\leq C_\beta\langle\xi\rangle^{-|\beta|}m(\xi)^{-1}.
\]
Hence, by Leibniz' rule and the comparison of $w$ and $m$,
\[
\left|\partial^\beta\left(\frac{w}{m}\right)(\xi)\right|
+
\left|\partial^\beta\left(\frac{m}{w}\right)(\xi)\right|
\leq C_\beta\langle\xi\rangle^{-|\beta|}.
\]
Then, the Mihlin--Hörmander theorem shows that both $w/m$
and $m/w$ are $L^p$-Fourier multipliers for
$1<p<\infty$, and consequently, for $u\in\mathcal S'(\mathbb R^n)$,
\[
\mathcal F^{-1}(w\widehat u)\in L^p(\mathbb R^n)
\quad\Longleftrightarrow\quad
\mathcal F^{-1}(m\widehat u)\in L^p(\mathbb R^n),
\]
and, whenever these conditions hold,
\[
\|\mathcal F^{-1}(w\widehat u)\|_{L^p}
\asymp
\|\mathcal F^{-1}(m\widehat u)\|_{L^p}.
\]
The Fourier characterization of
$H^{\rho_{\log}^s,p}(\mathbb R^n)$ and the definition of
$H_{\log}^{s,p,\lambda}(\mathbb R^n)$ now yield the
claimed equality with equivalence of norms. \qed
\section[The Gustavsson--Peetre method and related function spaces]{The Gustavsson--Peetre interpolation method and related function spaces}

We briefly recall the Gustavsson--Peetre $\pm$ method with a function parameter following
\cite[Section~6]{GustavssonPeetre}.
\begin{Defi}
    We define the class of functions  $\mathfrak{B}$ as all the continuous functions $\varphi:(0,\infty)\to (0,\infty)$ such that $\varphi(1)=1$ and $$\overline{\varphi}(\lambda):=\operatorname{sup}_{t>0}\frac{\varphi(\lambda t)}{\varphi(t)}<\infty,\,\lambda>0.$$ The function $\overline{\varphi}$ is called the \textit{dilation function} of $\varphi$.
\end{Defi}
\begin{Defi}[Gustavsson--Peetre $\pm$ method]\label{def:functor}
Let $(E_0,E_1)$ be a compatible Banach couple and let
$\varphi\in \mathfrak{B}$. For a sequence $\mathbf{x}=(x_m)_{m\in\Z}$ in
$E_0\cap E_1$, set
\[
 \mathfrak G_\varphi(\mathbf{x})=\max_{j=0,1}
 \sup_{\substack{F\subset\Z\ \mathrm{finite}\\|\varepsilon_m|\leq1}}
 \left\|\sum_{m\in F}\varepsilon_m
       \frac{2^{jm}}{\varphi(2^m)}x_m\right\|_{E_j}.
\]
The space $\langle E_0,E_1\rangle_\varphi$ consists of all
$x=\sum_{m\in\Z}x_m$, with convergence in $E_0+E_1$ and
$\mathfrak G_\varphi(\mathbf{x})<\infty$, equipped with the norm
\[
 \|x\|_{\langle E_0,E_1\rangle_\varphi}
 =\inf_{x=\sum x_m}\mathfrak G_\varphi(\mathbf{x}).
\]

\end{Defi}

For reflexive couples, the choice $\varphi(t)=t^\theta$ gives an equivalent method to Peetre's classical
spaces $\langle E_0,E_1\rangle_\theta$. See
\cite[Section~2]{Cwikel} for a modern exposition on these methods.

The original definition requires bounded finite signed sums, rather
than unconditional convergence of the weighted series. If both
endpoints are reflexive, these conditions coincide by
\cite[Lemma~2 and Theorem~5, pp.~159--160]{BessagaPelczynski}.
In that case finite truncations are dense in the sequence norm, so
the space also agrees with the completed-sequence formulation and
$E_0\cap E_1$ is dense in it.

In order for the method to be an interpolation functor we require the following additional hypotheses for functions in the class $\mathfrak{B}$:
For $\varphi\in \mathfrak{B}$, the functor $(E_0,E_1)\mapsto \langle E_0,E_1\rangle_{\varphi}$  preserves
boundedness of linear operators if \begin{equation}\label{eqn:condition}s_0(\overline{\varphi}):=\lim_{t\to 0^+}\frac{\log \overline{\varphi}(t)}{\log t}>0,\,s_1(\overline{\varphi}):=\lim_{t\to \infty}\frac{\log \overline{\varphi}(t)}{\log t}<1.\end{equation} These two quantities are called the \textit{Boyd indices} (see \cite{Boyd,CobosLibro}) of the dilation function. 
The following result from \cite[Proposition~6.1]{GustavssonPeetre} establishes the validity of the method as an interpolation functor.
\begin{Prop}[Gustavsson--Peetre]\label{prop:functor}
The construction is an exact interpolation functor of Banach spaces, i.e.,
for every linear operator $T:(E_0,E_1)\to(F_0,F_1)$ with positive endpoint
bounds $M_i$, one has $$\|T\|_{\langle E_0,E_1\rangle_\varphi\to \langle F_0,F_1\rangle_\varphi}\leq\max(M_0,M_1).$$ Moreover, 
\begin{equation}\label{eq:operator-bound}
 \|T\|_{\langle E_0,E_1\rangle_\varphi\to \langle F_0,F_1\rangle_\varphi}
 \lesssim M_0 \overline{\varphi}(M_1/M_0),\qquad
 \overline{\varphi}(\lambda):=\sup_{t>0}\frac{\varphi(\lambda t)}{\varphi(t)},\,\lambda>0.
\end{equation}
\end{Prop}
Another related interpolation method that will be useful in our development is the \textit{generalized real interpolation method} in the spirit of Gustavsson and Peetre (see \cite{Gustavsson} and the references therein). We briefly recall it and its main properties. Let $(E_0,E_1)$ be a compatible couple of Banach spaces, consider the sum $E_0+E_1$ with the norm $K(1,x)$, where $$K(t,x)=K(t,x;E_0,E_1):=\operatorname{inf}_{x=x_0+x_1}\left(\norm{x_0}_{E_0}+t\norm{x_1}_{E_1}\right),\,t>0.$$
\begin{Defi}[Generalized real interpolation method]
    Let $(E_0,E_1)$ be a compatible couple of Banach spaces. Let $\varphi\in \mathfrak{B}$ satisfying conditions $\ref{eqn:condition}$, $1\leq q\leq \infty$ and $u$ a measurable non-negative function of one variable. We define $$\Phi_{\varphi,q}(u):=\left(\int_0^\infty \left(\frac{u(t)}{\varphi(t)}\right)^q\frac{dt}{t}\right)^{1/q},$$ for $1\leq q<\infty$, and the corresponding one taking supremums for $q=\infty$. We define the spaces $(E_0,E_1)_{\varphi,q}$ as the set of elements $x\in E_0+E_1$ such that $$\Phi_{\varphi,q}(K(\cdot,x))<\infty.$$
\end{Defi} Note that for the choice $\varphi(t)=t^\theta,\,\theta\in (0,1)$, we recover the classical real interpolation method with parameter $\theta$. As it was proved in \cite[Theorem 2.1]{Gustavsson}, these methods are interpolation functors.
\begin{Teor}[Gustavsson]
    Let $T$ be a bounded linear operator from $E_j\to F_j$, with norm $M_j$, $j=0,1$. Then, $T:(E_0,E_1)_{\varphi,q}\to (F_0,F_1)_{\varphi,q}$ continuously with norm $$\norm{T}_{(E_0,E_1)_{\varphi,q}\to (F_0,F_1)_{\varphi,q}}\leq M_0\overline{\varphi}\left(\frac{M_1}{M_0}\right).$$
\end{Teor}
\subsection{Spaces of generalized smoothness}
\begin{Defi}[RO class]
    We define the class of functions \textbf{RO} as the Borel measurable functions $\alpha:[1,\infty)\to (0,\infty)$ such that there exist numbers $b>1$ and $c\geq 1$ (depending on $\alpha$), such that $$c^{-1}\leq \frac{\alpha(\lambda t)}{\alpha(t)}\leq c,\,t\geq 1,\,\lambda\in [1,b].$$
\end{Defi}
In order to introduce the Besov and Triebel-Lizorkin scales of generalized smoothness, we briefly recall some definitions:
\begin{Defi}
    A sequence of test functions $(\phi_j)_j\subset\mathcal{S}(\R^n)$ is called \textit{special} if, for some $N\in\N$, the Fourier transforms of its terms satisfy the following five properties: \begin{itemize}
        \item $\widehat{\phi_j}\geq 0$, for every $j$;
        \item $\operatorname{supp}\widehat{\phi_j}\subset \{\xi\in\R^n:2^{j-N}\leq |\xi|\leq 2^{j+N}\};$
        \item for every number $\delta\in (0,1)$ there exists $c_\delta$ such that, for each $j\in \N$, one has $\widehat{\phi_j}(\xi)\geq c_\delta$ when $(2-\delta)^{-N}2^j\leq |\xi|\leq (2-\delta)^N2^j$, and, in addition, $\widehat{\phi_0}(\xi)\geq c_\delta$ when $|\xi|\leq (2-\delta)^N$;
        \item there exists a number $c>0$ such that $c\leq \sum_{j=0}^\infty\widehat{\phi_j}(\xi)$, $\xi\in \R^n$;
        \item for every multi-index $\beta$, there exists a number $c_\beta>0$ such that $2^{j|\beta|}|\partial ^\beta\widehat{\phi_j}(\xi)|\leq c_\beta$, for all $\xi\in\R^n$ and $j\in \N$.
   \end{itemize}
\end{Defi}
If we choose any $N\in\N$ and a function $\eta\in \mathcal{S}(\R^n)$ such that $\operatorname{supp}\eta \subset \{\xi\in \R^n: 2^{-N}\leq |\xi|\leq 2^N\}$ and $\eta(\xi)>0$ in $2^{-N}\leq |\xi|\leq 2^N$, if we set $\widehat{\eta_j}(\xi):=\eta(2^{-j}\xi)$ for $j>0$ and some arbitrary function $\widehat{\eta}_0\in \mathcal{S}(\R^n)$ such that $\operatorname{supp}\widehat{\eta_0}\subset\{\xi\in \R^n:|\xi|\leq 2^N\}$ and $\widehat{\eta_0}(\xi)>0$ whenever $|\xi|<2^N$, taking inverse Fourier transforms yields a special sequence.

We fix a special sequence $(\phi_j)_j$.
\begin{Defi}[Besov and Triebel-Lizorkin spaces with generalized smoothness]Let $\alpha\in \textbf{RO}$ and $p,q\in[1,\infty]$. We define the Besov spaces with generalized smoothness $B^{\alpha(\cdot)}_{p,q}(\R^n)$ as the linear space of distributions $u\in \mathcal{S}'(\R^n)$ satisfying $$\norm{u}_{B^{\alpha(\cdot)}_{p,q}(\R^n)}:=\left(\sum_{j=0}^\infty\alpha^q(2^j)\norm{\phi_j*u}_p^q\right)^{1/q}<\infty,$$ for $q<\infty$, and the corresponding definition when $q=\infty$. The Triebel-Lizorkin spaces with generalized smoothness $F^{\alpha(\cdot)}_{p,q}(\R^n)$, for $1\leq p,q<\infty$, are defined as the linear space of distributions $u\in \mathcal{S}'(\R^n)$ satisfying
    $$\norm{u}_{F^{\alpha(\cdot)}_{p,q}(\R^n)}:=\left(\int_{\R^n}\left(\sum_{j=0}^\infty \alpha^q(2^j)|\phi_j*u(x)|^q\right)^{p/q}\,dx\right)^{1/p}<\infty.$$ Neither definition depends on the choice of the special sequence. Furthermore, $F^{\alpha(\cdot)}_{p,p}(\R^n)=B^{\alpha(\cdot)}_{p,p}(\R^n)$, whenever $1\leq p<\infty$, and for the choice $\alpha(t)=t^s$, $s\in \R$, the spaces become the classical Besov and Triebel-Lizorkin spaces.
\end{Defi}
The following result from \cite[Theorem 4.1]{CobosLibro} recalls the main continuous embeddings on these scales:
\begin{Teor}\label{EmbeddingsNL}
    Let $\alpha\in \textbf{RO}$. Then, the following continuous embeddings hold:
    \begin{itemize}
        \item $F^{\alpha(\cdot)}_{p,q_0}(\R^n)\xhookrightarrow{}F^{\alpha(\cdot)}_{p,q_1}(\R^n)$, $1<p<\infty, 1<q_0\leq q_1<\infty$;
        \item $B^{\alpha(\cdot)}_{p,q_0}(\R^n)\xhookrightarrow{}B^{\alpha(\cdot)}_{p,q_1}(\R^n)$, $1\leq p<\infty, 1<q_0\leq q_1\leq \infty$;
        \item $B^{\alpha(\cdot)}_{p,q}(\R^n)\xhookrightarrow{}F^{\alpha(\cdot)}_{p,q}(\R^n)$, $1<q\leq p<\infty$;
        \item $F^{\alpha(\cdot)}_{p,q}(\R^n)\xhookrightarrow{}B^{\alpha(\cdot)}_{p,q}(\R^n)$, $1<p\leq q<\infty$;
        \item $B^{\alpha(\cdot)}_{p,1}(\R^n)\xhookrightarrow{}F^{\alpha(\cdot)}_{p,q}(\R^n)\xhookrightarrow{}B^{\alpha(\cdot)}_{p,\infty}(\R^n)$, $1<p,q<\infty$. 
        \end{itemize}
        Moreover, if $\alpha_0,\alpha_1\in \textbf{RO}$ are such that $\alpha_1/\alpha_0\in L^1_*(1,\infty)$, we have that \begin{itemize}
            \item $B^{\alpha_0(\cdot)}_{p,\infty}(\R^n)\xhookrightarrow{}B^{\alpha_1(\cdot)}_{p,1}(\R^n)$, $1\leq p\leq \infty$;
            \item $F^{\alpha_0(\cdot)}_{p,q_0}(\R^n)\xhookrightarrow{}F^{\alpha_1(\cdot)}_{p,q_1}(\R^n)$, $1<p,q_0,q_1<\infty$;
        \end{itemize}
    where $\alpha_1/\alpha_0\in L_*^1(1,\infty)$ means that $\int_1^\infty \frac{\alpha_1(t)}{\alpha_0(t)}\frac{dt}{t}<\infty$.
\end{Teor}
If we consider a broader class of \textbf{RO} functions, we can define Sobolev spaces of generalized smoothness analogously as in the fractional case with Bessel potential spaces.
\begin{Defi}[Sobolev spaces with generalized smoothness]
    Let $\alpha\in \textbf{RO}$ such that $\alpha\in C^\infty([1,\infty))$ and $$|\alpha^{m)}(t)|\lesssim t^{-m}\alpha(t),\,t\geq 1,\,m\in \N.$$ Then, for every $1<p<\infty$, we define the \textit{Sobolev spaces with generalized smoothness} as $$H^{\alpha(\cdot),p}(\R^n):=\{u\in \mathcal{S}'(\R^n): \mathfrak{F}^{-1}\left(\alpha(\langle \xi\rangle)\widehat{u}(\xi)\right)\in L^p(\R^n)\},$$ with the norm $$\norm{u}_{H^{\alpha(\cdot),p}}:=\norm{\mathfrak{F}^{-1}\left(\alpha(\langle \xi\rangle)\widehat{u}(\xi)\right)}_p.$$ 
\end{Defi}
Actually, the generalized smoothness Bessel scale is well defined if the function $\alpha(\langle\xi\rangle)$ is a
pointwise multiplier on $\mathcal{S}'(\R^n)$, which is more general than the previous hypothesis.
A function $\alpha\in \textbf{RO}$ under the previous hypothesis will be said to be a multiplier \textbf{RO}-function or, alternatively, to belong to the \textbf{MRO} class.
Under those hypotheses, the $H^{\alpha(\cdot),p}$-spaces coincide with equivalence of the norms with the generalized smoothness Triebel-Lizorkin spaces $F^{\alpha(\cdot)}_{p,2}(\R^n)$ (see \cite[Theorem 1]{AnopMurach}). This fact follows as in the classical theory by means of the operator $$J^{\alpha(\cdot)}:f\mapsto \mathfrak{F}^{-1}\left(\alpha(\langle\xi\rangle)\widehat{f}(\xi)\right),\,f\in \mathcal{S}'(\R^n),$$ which is a lifting from $H^{\alpha(\cdot),p}(\R^n)\to L^p(\R^n)$ \cite[Theorem 3.4]{CobosLibro}. Moreover, by \cite[Theorem 3.1]{CobosLibro}, for $\alpha_0\in \textbf{RO}, \alpha_1\in \textbf{MRO}$, $J^{\alpha_1}$ is an isomorphism \begin{align}\label{lifting}
    F^{\alpha_0(\cdot)}_{p,q}(\R^n)&\to F^{(\alpha_0/\alpha_1)(\cdot)}_{p,q}(\R^n),\,1<p,q<\infty,\\
    B^{\alpha_0(\cdot)}_{p,q}(\R^n)&\to B^{(\alpha_0/\alpha_1)(\cdot)}_{p,q}(\R^n),\,1\leq p,q\leq \infty.
\end{align}
Note that for the choice $\alpha_s(t)=t^s$, $s\in \R$, $\alpha_s\in \textbf{RO}$ and $H^{\alpha_s(\cdot),p}(\R^n)=H^{s,p}(\R^n)$, the classical Bessel potential spaces. 
Generalized smoothness Besov and Triebel-Lizorkin spaces behave as expected under the complex interpolation method. In fact, we have the following result \cite[Theorem 2]{AnopMurach}:
\begin{Teor}\label{complex}
    Let $\theta\in (0,1), p_0,p_1,q_0,q_1\in [1,\infty)$ and $\alpha_0,\alpha_1\in \textbf{RO}.$ We define the function $$\alpha_\theta(t):=\alpha_0(t)^{1-\theta}\alpha_1(t)^\theta,\,t\geq 1.$$ Then, we have that $\alpha_\theta\in \textbf{RO}$, and $$[F^{\alpha_0(\cdot)}_{p_0,q_0}(\R^n),F^{\alpha_1(\cdot)}_{p_1,q_1}(\R^n)]_\theta=F^{\alpha_\theta(\cdot)}_{p_\theta,q_\theta}(\R^n),\,[B^{\alpha_0(\cdot)}_{p_0,q_0}(\R^n),B^{\alpha_1(\cdot)}_{p_1,q_1}(\R^n)]_\theta=B^{\alpha_\theta(\cdot)}_{p_\theta,q_\theta}(\R^n),$$ where $$\frac{1}{p_\theta}=\frac{1-\theta}{p_0}+\frac{\theta}{p_1},\,\frac{1}{q_\theta}=\frac{1-\theta}{q_0}+\frac{\theta}{q_1}$$ with equivalence of norms. If $\alpha_0,\alpha_1\in\textbf{MRO}$, the analogous result holds for the generalized smoothness Bessel spaces.
\end{Teor}

An important observation in our development is that for every function $\alpha\in \textbf{RO}$, there exist real numbers $s_0<s_1$ such that \begin{equation}\label{eqn:cociente}\lambda^{s_0}\lesssim \frac{\alpha(\lambda t)}{\alpha(t)}\lesssim \lambda^{s_1},\,t,\lambda\in [1,\infty),\end{equation} as it was proved in \cite[Appendix 1, Theorem 2]{Seneta}. The parameters \begin{align*}
    \sigma_0(\alpha)&:=\operatorname{sup}\{s_0\in \R: \text{the left-hand inequality in \ref{eqn:cociente} holds}\},\\
    \sigma_1(\alpha)&:=\operatorname{inf}\{s_1\in \R: \text{the right-hand inequality in \ref{eqn:cociente} holds}\},
\end{align*} are called \textit{the lower and upper Matuszewska indices} (see \cite{japo}) of the function $\alpha$.

Classes $\mathfrak{B}$ and \textbf{RO} are deeply connected. First of all, we have that the scales of generalized smoothness Besov and Triebel-Lizorkin spaces $B^{\alpha(\cdot)}_{p,q}, F^{\alpha(\cdot)}_{p,q}$ are the same with equivalence of norms if we consider $\alpha\in \mathfrak{B}$ or $\alpha\in \textbf{RO}$. The key here is that from a continuous function $\alpha\in \textbf{RO}$, we can define a function $\alpha^*\in \mathfrak{B}$ by $$\alpha^*(t):=\begin{cases}
    \alpha(t)/\alpha(1),\,t\geq 1,\\
    \alpha(1)/\alpha(t^{-1}),\,0<t<1.
\end{cases}.$$ For this function, we have that $$\sigma_0(\alpha)=s_0(\overline{\alpha^*}),\,\sigma_1(\alpha)=s_1(\overline{\alpha^*}),$$ as noted in \cite[(65)]{AnopMurach}.

Moreover, we have the following interpolation result for the Gustavsson-Peetre method in connection with \textbf{MRO} functions \cite[Theorem 16]{AnopMurach}
\begin{Teor}
    Let $\alpha_0,\alpha_1\in \textbf{RO}$, $1\leq p\leq \infty$, $1<q\leq \infty$ and $\varphi\in \mathfrak{B}$ satisfying \ref{eqn:condition}. If we define $$\alpha_\varphi(t):=\frac{\alpha_0(t)}{\varphi(\alpha_0(t)/\alpha_1(t))},\,t\geq 1,$$ then $\alpha_\varphi\in \textbf{RO}$ and $$\langle F^{\alpha_0(\cdot)}_{p,q}(\R^n),F^{\alpha_1(\cdot)}_{p,q}(\R^n)\rangle_\varphi=F^{\alpha_\varphi(\cdot)}_{p,q}(\R^n),$$ with equivalence of norms.
\end{Teor}
Furthermore, choosing $q=2$, we can construct every space on the $H^{\alpha(\cdot),p}$-scale as an interpolation space between purely fractional Bessel potential spaces by means of the Gustavsson-Peetre method \cite[Theorem 5, Theorem 17]{AnopMurach}:
\begin{Teor}\label{TEOREMA}
    Let $\alpha\in \textbf{MRO}$ and $1<p<\infty$. For arbitrary real numbers $s_0<\sigma_0(\alpha)$ and $s_1>\sigma_1(\alpha)$, we define $$\varphi_\alpha(t):=\begin{cases}
        t^{-s_0/(s_1-s_0)}\widetilde{\alpha}(t^{1/(s_1-s_0)}),\,t\geq 1;\\
        1/\varphi_\alpha(1/t),\,0<t<1.
    \end{cases},$$ where $\widetilde{\alpha}$ is a continuous \textbf{MRO} function such that $\widetilde{\alpha}(1)=1$ and is equivalent to $\alpha$ for $t\geq 1$ (it is always possible to choose such a function). Then, $\varphi_\alpha\in \mathfrak{B}$ satisfying \ref{eqn:condition}, and $$H^{\alpha(\cdot),p}(\R^n)=\langle H^{s_0,p}(\R^n),H^{s_1,p}(\R^n)\rangle_{\varphi_\alpha},$$ with equivalence of norms.
\end{Teor}
Note that the two previous interpolation results are originally stated in full generality on the Triebel-Lizorkin scales. If we now look at the Besov scale, we also have analogous interpolation results by the generalized real method. The first one is \cite[Theorem 4]{AnopMurach}
\begin{Teor}
    Let $q\in [1,\infty],$ $\alpha_0,\alpha_1\in \textbf{RO}$ with $\sigma_0(\alpha_0/\alpha_1)>0$ or $\sigma_1(\alpha_0/\alpha_1)<0,$ and $\varphi\in \mathfrak{B}$ satisfying \ref{eqn:condition}. If we define $$\alpha_\varphi(t):=\frac{\alpha_0(t)}{\varphi(\alpha_0(t)/\alpha_1(t))},\,t\geq 1,$$ then $\alpha_\varphi\in \textbf{RO}$ and \begin{align*}\left(F^{\alpha_0(\cdot)}_{p,q_0}(\R^n),F^{\alpha_1(\cdot)}_{p,q_1}(\R^n)\right)_{\varphi,q}&=B^{\alpha_\varphi(\cdot)}_{p,q}(\R^n),\,p,q_0,q_1\in [1,\infty),\\
    \left(B^{\alpha_0(\cdot)}_{p,q_0}(\R^n),B^{\alpha_1(\cdot)}_{p,q_1}(\R^n)\right)_{\varphi,q}&=B^{\alpha_\varphi(\cdot)}_{p,q}(\R^n),\,p,q_0,q_1\in [1,\infty],
    \end{align*} with equivalence of norms. Furthermore, if $\alpha_0,\alpha_1\in \textbf{MRO}$, then $\alpha_\varphi\in \textbf{MRO}$ and $$\left(H^{\alpha_0(\cdot),p}(\R^n),H^{\alpha_1(\cdot),p}(\R^n)\right)_{\varphi,q}=B^{\alpha_\varphi(\cdot)}_{p,q}(\R^n),\,1<p<\infty,$$ with equivalence of norms.
\end{Teor}
We also have the analogue of Theorem \ref{main} for the Besov scale \cite[Theorem 5]{AnopMurach}:
\begin{Teor}\label{TEOREMA2}
     Let $\alpha\in \textbf{RO}$ and $1\leq p,q\leq\infty$. For arbitrary real numbers $s_0<\sigma_0(\alpha)$ and $s_1>\sigma_1(\alpha)$, we define $$\varphi_\alpha(t):=\begin{cases}
        t^{-s_0/(s_1-s_0)}\widetilde{\alpha}(t^{1/(s_1-s_0)}),\,t\geq 1;\\
        1/\varphi_\alpha(1/t),\,0<t<1.
    \end{cases},$$ where $\widetilde{\alpha}$ is a continuous \textbf{RO} function such that $\widetilde{\alpha}(1)=1$ and is equivalent to $\alpha$ for $t\geq 1$ (it is always possible to choose such a function). Then, $\varphi_\alpha\in \mathfrak{B}$ satisfying \ref{eqn:condition}, and $$B^{\alpha(\cdot)}_{p,q}(\R^n)=\left(B^{s_0}_{p,q_0}(\R^n),B^{s_1}_{p,q_1}(\R^n)\right)_{\varphi_\alpha,q},\,q_0,q_1\in [1,\infty],$$ with equivalence of norms.
\end{Teor}

\section{Nonlocal Sobolev spaces as interpolation spaces}
In this section we study the interpolation of nonlocal Sobolev spaces for general kernels. First of all we will identify them as an interpolation space by the $\pm$-method. Then, we will study the interpolation between those spaces for different kernels and finally we will introduce the bounded domain version via interpolation.
\subsection{Interpolation by the Gustavsson-Peetre method}
Our main result reads as follows:
\begin{Teor}\label{main}
    Let $\rho$ satisfy $(H0)-(H4)$ with compact support and $1<p<\infty$. We define $$\alpha_\rho(t):=w_\rho(\sqrt{t^2-1}e_1),\,t>1,\,\alpha_\rho(1)=1,$$ and consider $$\varphi_\rho(t):=\begin{cases}
        \alpha_\rho(t^{-1})^{-1},\,0<t<1,\\
        \alpha_\rho(t),\,t\geq 1.
    \end{cases}$$ Then, $\alpha_\rho\in \textbf{MRO}$ and $$0<\sigma\leq\sigma_0(\alpha_\rho)\leq\sigma_1(\alpha_\rho)\leq \gamma<1,$$ and the function $\varphi_\rho$ belongs to $\mathfrak{B}$ and satisfies \ref{eqn:condition}. Therefore $$H^{\rho,p}(\R^n)=\langle L^p(\R^n),W^{1,p}(\R^n)\rangle_{\varphi_\rho}.$$
\end{Teor}
\noindent{}\textbf{Proof:} Suppose that the first part of the theorem is true. Hence, by virtue of Theorem \ref{TEOREMA}, for the admissible choices $s_0=0$ and $s_1=1$, we have that with equivalence of norms, $$\langle L^p(\R^n),W^{1,p}(\R^n)\rangle_{\varphi_\rho}=H^{\alpha_\rho(\cdot),p}(\R^n)=\{u\in \mathcal{S}'(\R^n): \mathfrak{F}^{-1}\left(\alpha_\rho(\langle\xi\rangle)\widehat{u}(\xi)\right)\in L^p(\R^n)\}.$$ Now, observe that by radiality of $w_\rho$, $$\alpha_\rho(\langle \xi\rangle)=w_\rho(e_1\sqrt{\langle\xi\rangle^2-1})=w_\rho(e_1|\xi|)=w_\rho(\xi),$$ and since $\langle\xi\rangle\geq 1$ for every $\xi\in \R^n$, this identification holds for every $\xi\in \R^n$. Moreover, by the potential characterization of $H^{\rho,p}(\R^n)$ in terms of $w_\rho$, we have that $H^{\rho,p}(\R^n)=H^{\alpha_\rho(\cdot),p}(\R^n)$, as we want to prove. Hence, in order to establish the result we only need to prove that the functions $\alpha_\rho,\varphi_\rho$ satisfy the hypothesis of Theorem \ref{TEOREMA}.

Smoothness of $\alpha_\rho$ for $t>1$ follows directly from the previous lemma. At $t=1$, the compact support of $Q_\rho$ allows us to extend the function analytically $\alpha_\rho$ to $\alpha_\rho(1)=1$. Therefore $\alpha_\rho\in C^\infty([1,\infty))$. 

To control the derivatives of $\alpha_\rho$, we first note that $\alpha_\rho(t)=w_\rho(g(t)e_1)$ where $g(t)=\sqrt{t^2-1}$. For $t>1$, we can apply Fa\`a di Bruno's formula to get that for every $k\in \N$, \begin{align*}
    \alpha_\rho^{k)}(t)=\sum_{\substack{m_1,m_2,\ldots,m_k\geq 0,\\m_1+2m_2+\cdots+km_k=k}}\frac{k!}{m_1!\cdot\dots\cdot m_k!}\left(\partial_1^{m_1+\ldots+m_k}w_\rho\right)(g(t)e_1)\prod_{j=1}^k\left(\frac{g^{j)}(t)}{j!}\right)^{m_j}.
\end{align*} By the previous lemma, \begin{align*}|(\partial_1^{m_1+\ldots+m_k}w_\rho)(g(t)e_1)|&\leq C_{m_1+\ldots+m_k}\langle g(t)e_1\rangle^{-(m_1+\ldots+m_k)}w_\rho(g(t)e_1)\\
&=C_{m_1+\ldots+m_k}t^{-(m_1+\ldots+m_k)}\alpha_\rho(t).\end{align*} For $t\geq 2$, we can find for each $j=1,\ldots,k$, a positive constant $C_j$ such that $|g^{j)}(t)|\leq C_j t^{1-j}$, and hence $$\prod_{j=1}^k \left(\frac{|g^{j)}(t)|}{j!}\right)^{m_j}\leq \prod_{j=1}^k \left(\frac{C_j}{j!}\right)^{m_j}(t^{1-j})^{m_j}\leq C_k t^{\sum_{j=1}^k(1-j)m_j}=C_kt^{(m_1+\ldots+m_k)-k}.$$ Putting both estimates together we obtain that  \begin{align*}
    |\alpha_\rho^{k)}(t)|&\leq \sum_{\substack{m_1,m_2,\ldots,m_k\geq 0,\\m_1+2m_2+\cdots+km_k=k}}\frac{k!}{m_1!\cdot\dots\cdot m_k!}C_{m_1+\ldots+m_k}t^{-(m_1+\ldots+m_k)}|\alpha_\rho(t)|C_kt^{(m_1+\ldots+m_k)-k}\\
    &=\widetilde{C}_kt^{-k}\alpha_\rho(t),\,t\geq 2,
\end{align*} for some positive constant $\widetilde{C}_k$. On the interval $t\in [1,2]$, since $\alpha_\rho\in C^\infty([1,\infty))$, there exists $M_k:=\operatorname{max}_{t\in [1,2]}|\alpha_\rho^{(k)}(t)|<\infty.$ Moreover, $t^k\leq 2^k$ and $\alpha_\rho(t)\geq 1$, hence $$|\alpha_\rho^{k)}(t)|\leq M_k2^kt^{-k}\alpha_\rho(t),\,t\in [1,2].$$ Finally, $$|\alpha_\rho^{k)}(t)|\leq \operatorname{max}\{\widetilde{C}_k,M_k2^k\}t^{-k}\alpha_\rho(t),\,t\geq 1.$$

Now, since $|\alpha_\rho'(t)|\leq C t^{-1}\alpha_\rho(t),\,t\geq 1$, and $\alpha_\rho(t)>0$, we have that $$-Ct^{-1}\leq \frac{\alpha_\rho'(t)}{\alpha_\rho(t)}\leq Ct^{-1},$$ and integrating both sides between $t$ and $\lambda t$, we get that $$-C\log \lambda\leq \log \frac{\alpha_\rho(\lambda t)}{\alpha_\rho(t)}\leq C\log \lambda,$$ which implies that $$\lambda^{-C}\leq \frac{\alpha_\rho(\lambda t)}{\alpha_\rho(t)}\leq \lambda^C.$$
If we restrict ourselves to $1\leq \lambda\leq 2$, we get that $$2^{-C}\leq \frac{\alpha_\rho(\lambda t)}{\alpha_\rho(t)}\leq 2^C,$$ i.e., $\alpha_\rho\in \textbf{RO}$ for $b=2$ and $c=2^C$. Furthermore, by the previous estimates, $\alpha_\rho\in \textbf{MRO}$ with $\alpha_\rho(1)=1$.

The next step is to relate Matuszewska indices with the parameters $\sigma,\gamma$ from hypothesis $(H3)-(H4)$. Let $a(t):=t^{1-n}\overline{\rho}(1/t)$. Since $\alpha_\rho(t)=w_\rho(g(t)e_1)$ and $w_\rho(\xi)\asymp |\xi|^{1-n}\overline{\rho}(1/|\xi|)$ for $|\xi|$ large enough, we get that, for sufficiently large $t$, $\alpha_\rho(t)\asymp a(g(t))$. Let us fix $t_0:=\operatorname{min}\{1,\varepsilon\}$. By the hypothesis $(H3)$, we have that for every $0<s\leq t<t_0$, $$s^{n+\sigma-1}\overline{\rho}(s)\geq c t^{n+\sigma-1}\overline{\rho}(t)>0,$$ which implies that $$\frac{\overline{\rho}(s)}{\overline{\rho}(t)}\geq c\left(\frac{t}{s}\right)^{n+\sigma-1}.$$ By analogous reasoning, using hypothesis $(H4)$ gives that $$\frac{\overline{\rho}(s)}{\overline{\rho}(t)}\leq C\left(\frac{t}{s}\right)^{n+\gamma-1}.$$ Now, we take $r=t^{-1}\geq 1$ and $s=(\lambda r)^{-1}$ for some $\lambda\geq 1.$ Then, $s\leq t\leq t_0$, $t/s=\lambda$, and hence $$c\lambda^{n+\sigma-1}\leq \frac{\overline{\rho}(1/(\lambda r))}{\overline{\rho}(1/r)}\leq C\lambda^{n+\gamma-1}.$$ Multiplying each side by $\lambda^{1-n}$ we get that $$c\lambda^{\sigma}\leq \lambda^{1-n}\frac{\overline{\rho}(1/(\lambda r))}{\overline{\rho}(1/r)}\leq C\lambda^{\gamma},$$ with $$\lambda^{1-n}\frac{\overline{\rho}(1/(\lambda r))}{\overline{\rho}(1/r)}=\frac{a(\lambda r)}{a(r)}.$$ Now, for $t\geq 2$, it is clear that $1\leq t/g(t)\leq 2$. By the previous estimate for $r=g(t)$ and $\lambda=t/g(t)$, we get that 
$$c\leq c\left(\frac{t}{g(t)}\right)^\sigma\leq \frac{a(t)}{a(g(t))}\leq C\left(\frac{t}{g(t)}\right)^\gamma\leq C2^\gamma,$$ which means that $a(t)\asymp a(g(t))$, and hence $\alpha_\rho(t)\asymp a(t)$ for $t\geq 2$ large enough, so we obtain that $$c\lambda^{\sigma}\leq \frac{\alpha_\rho(\lambda t)}{\alpha_\rho(t)}\leq C\lambda^{\gamma},$$ for $t$ large enough. Let us say that this estimate holds for $t\geq T$ for some $T\geq 2$. Now, if $1\leq t\leq \lambda t\leq T$, both values of $\alpha_\rho$ lie in a fixed
positive compact range and $\lambda\leq T$, so the same bounds hold after enlarging the constants. Finally, if $1\leq t\leq T\leq\lambda t$, we write $$\frac{\alpha_\rho(\lambda t)}{\alpha_\rho(t)}=\frac{\alpha_\rho(\lambda t)}{\alpha_\rho (T)}\frac{\alpha_\rho(T)}{\alpha_\rho(t)}.$$ For the first quotient, we have by the previous estimates that $$c\left(\frac{\lambda t}{T}\right)^\sigma\leq \frac{\alpha_\rho(\lambda t)}{\alpha_\rho (T)}\leq C\left(\frac{\lambda t}{T}\right)^\gamma,$$ and for the second quotient $$\frac{\alpha_\rho(T)}{M}\leq \frac{\alpha_\rho(T)}{\alpha_\rho(t)}\leq \frac{\alpha_\rho(T)}{m},$$ where $m:=\operatorname{min}_{[1,T]}\alpha_\rho$ and $M:=\operatorname{max}_{[1,T]}\alpha_\rho$. Then, it follows that $$\frac{c\alpha_\rho(T)}{MT^\sigma}\lambda^\sigma\leq \frac{\alpha_\rho(\lambda t)}{\alpha_\rho(t)}\leq \frac{C\alpha_\rho(T)}{mT^\gamma}\lambda^\gamma,$$ so we can extend the bounds for the quotient $\alpha_\rho(\lambda t)/\alpha_\rho(t)$ for every $t\geq 1$ and $\lambda\geq 1$. From there, it is straightforward to see that $$0<\sigma\leq \sigma_0(\alpha_\rho)\leq \sigma_1(\alpha_\rho)\leq \gamma<1,$$ as we wanted to prove. In particular this proves that the particular choices $s_0=0$ and $s_1=1$ are in the hypothesis of Theorem \ref{TEOREMA}, and hence that the couple $(H^{0,p}(\R^n),H^{1,p}(\R^n))=(L^p(\R^n),W^{1,p}(\R^n))$ is admissible for interpolation by the Gustavsson-Peetre method. 

To conclude the result it only remains to prove that the function $\varphi_\rho$ defined in the statement of the theorem belongs to $\mathfrak{B}$ and satisfies \ref{eqn:condition}. Clearly, $\varphi_\rho$ is a continuous and positive function on $(0,\infty)$ such that $\varphi_\rho(1)=1$. An argument analogous to the preceding ones gives us that $$c\lambda^\sigma\leq \frac{\varphi_\rho(\lambda t)}{\varphi_\rho(t)}\leq C\lambda^\gamma,\,t>0,\,\lambda\geq 1.$$ If $\lambda\geq 1$, then it follows that $\overline{\varphi_\rho}(\lambda)\leq C\lambda^\gamma$. If $\lambda<1$, then $\lambda^{-1}>1$ and $$c\lambda^{-\sigma}\leq \frac{\varphi_\rho (\lambda^{-1}\lambda t)}{\varphi_\rho(\lambda t)}\leq C\lambda^{-\gamma}\implies \frac{\varphi_\rho (\lambda t)}{\varphi_\rho(t)}\leq c^{-1}\lambda^{\sigma},$$ thus $\overline{\varphi_\rho}(\lambda)\leq c^{-1}\lambda^\sigma.$ Then, we have the estimate $$\overline{\varphi_\rho}(\lambda)\leq \begin{cases}
    c^{-1}\lambda^\sigma,\,0<\lambda<1,\\
    C\lambda^\gamma,\,\lambda\geq 1
\end{cases},$$ which implies that $\varphi_\rho\in \mathfrak{B}$. Furthermore, from the last estimate it is straightforward to see that $$s_0(\overline{\varphi_\rho})=\lim_{\lambda\to 0^+}\frac{\log \overline{\varphi_\rho}(\lambda)}{\log\lambda}\geq \lim_{\lambda\to 0^+}\frac{\log c^{-1}+\sigma\log\lambda}{\log\lambda}=\sigma>0,$$ (where we have reversed the inequality since $\log \lambda<0$ for $\lambda$ close to $0$ enough), and $$s_1(\overline{\varphi_\rho})=\lim_{\lambda\to \infty}\frac{\log \overline{\varphi_\rho}(\lambda)}{\log\lambda}\leq \lim_{\lambda\to \infty}\frac{\log C+\gamma\log\lambda}{\log\lambda}=\gamma<1,$$ and hence $\varphi_\rho$ satisfies \ref{eqn:condition}. Then, the result follows. \qed

As a direct consequence of this identification we can compute the dual space of $H^{\rho,p}(\R^n)$, which agrees with the result obtained in \cite[Theorem 3.6]{Happ2026}. The result follows from \cite[Theorem 3.1.10]{Farkas}, with the identification of the spaces $F^{\sigma,N}_{p,q}$ introduced there with $F^{\alpha(\cdot)}_{p,q}$ for the admissible choices $N_j=2^j$ and $\sigma_j=\alpha(2^j)$, $j\in \N$. See the discussion in \cite[p. 8]{AnopMurach}.
Then, the result reads as follows:
\begin{Coro}
    Let $\rho$ be a radial kernel with compact support satisfying $(H0)-(H4)$. Then, $$\left(H^{\rho,p}(\R^n)\right)^*=H^{1/\rho,p'}(\R^n):=\{u\in \mathcal{S}'(\R^n): \mathfrak{F}^{-1}\left(w_\rho(\xi)^{-1}\widehat{u}(\xi)\right)\in L^{p'}(\R^n)\},\,1<p<\infty,$$ where $1/p+1/p'=1.$
\end{Coro}
\noindent{}\textbf{Proof:}
By Theorem \ref{TEOREMA}, we have the identification $H^{\rho,p}(\R^n)=H^{\alpha_\rho(\cdot),p}(\R^n)=F^{\alpha_\rho(\cdot)}_{p,2}(\R^n).$ Now, by \cite[Theorem 3.1.10]{Farkas}, we have that $$\left(F^{\alpha_\rho(\cdot)}_{p,2}(\R^n)\right)^*=F^{1/\alpha_\rho(\cdot)}_{p',2}=\{u\in \mathcal{S}'(\R^n): \mathfrak{F}^{-1}\left(\alpha_\rho(\langle\xi\rangle)^{-1}\widehat{u}(\xi)\right)\in L^{p'}(\R^n)\}.$$ Since $$\frac{1}{\alpha_\rho(\langle\xi\rangle)}=\frac{1}{w_\rho(\xi)},$$ we have that the dual of $H^{\rho,p}(\R^n)$ is given by $$H^{1/\rho,p'}(\R^n):=F^{1/\alpha_\rho(\cdot)}_{p',2}=\{u\in \mathcal{S}'(\R^n): \mathfrak{F}^{-1}\left(w_\rho(\xi)^{-1}\widehat{u}(\xi)\right)\in L^{p'}(\R^n)\},$$ as we wanted to prove.

\subsection{Complex interpolation of general nonlocal spaces}
In this section we give a result of interpolation between nonlocal Sobolev spaces by Calder\'on's complex method. Our result is an adaptation of the more abstract setting of Theorem \ref{complex} from \cite[Theorem 2]{AnopMurach}. This theorem states sort of a particular result on reiteration of the Gustavsson-Peetre method with Calder\'on's complex one. We first verify that the class of kernels under consideration is
stable under geometric interpolation.

\begin{lema}
\label{lem:geometric-kernel-stability}
Let $\rho_0,\rho_1$ be compactly supported radial kernels satisfying
\textup{(H0)--(H4)}, and let $0<\theta<1$. Then
\[
    \rho_\theta:=\rho_0^{1-\theta}\rho_1^\theta
\]
is a compactly supported radial kernel satisfying
\textup{(H0)--(H4)}. More precisely, if $\sigma_j,\gamma_j$ are
admissible exponents in \textup{(H3)--(H4)} for $\rho_j$, $j=0,1,$ then
\[
    \sigma_\theta=(1-\theta)\sigma_0+\theta\sigma_1,
    \qquad
    \gamma_\theta=(1-\theta)\gamma_0+\theta\gamma_1
\]
are admissible exponents for $\rho_\theta$.
\end{lema}
\noindent{}\textbf{Proof:}
Write $f_j(r)=r^{n-2}\overline\rho_j(r)$, $j=0,1$. Then
\[
    f_\theta(r):=r^{n-2}\overline\rho_\theta(r)
    =f_0(r)^{1-\theta}f_1(r)^\theta.
\]
Radiality, nonnegativity and compact support are immediate.
H\"older's inequality with respect to the measure
$d\nu(x)=\min\{1,|x|^{-1}\}\,dx$ gives
\[
    \int_{\mathbb R^n}\rho_\theta\,d\nu
    \leq
    \left(\int_{\mathbb R^n}\rho_0\,d\nu\right)^{1-\theta}
    \left(\int_{\mathbb R^n}\rho_1\,d\nu\right)^\theta
    <\infty.
\]
Both kernels are bounded below by a positive constant on a common
sufficiently small ball, so the same holds for $\rho_\theta$.
This proves \textup{(H0)}.

Since $f_0,f_1$ are nonincreasing, so is $f_\theta$. Moreover,
for sufficiently small $r>0$,
\[
    f_\theta(r)
    \leq \mu_0^{1-\theta}\mu_1^\theta f_\theta(r/2),
\]
where $0<\mu_j<1$ are the constants in \textup{(H1)}.
Thus \textup{(H1)} holds with
$\mu_\theta=\mu_0^{1-\theta}\mu_1^\theta<1$.

We next justify smoothness at the endpoints of the supports.
We claim that if $f\in C^\infty((0,\infty))$ is nonnegative,
nonincreasing and vanishes for all sufficiently large $r$, then
$f^a\in C^\infty((0,\infty))$ for every $a>0$.
The assertion is immediate if $f\equiv0$. Otherwise, let
$R=\sup\{r>0:f(r)>0\}$, so that $f>0$ on $(0,R)$ and
$f=0$ on $[R,\infty)$. Fix $k\geq1$ and an integer
$m>\max\{k,k/a\}$. Taylor's formula at $m$ distinct, equally
spaced points of $[r,r+h]$, followed by inversion of the resulting
Vandermonde matrix, yields
\[
    |f^{(\ell)}(r)|
    \leq C_m\bigl(h^{-\ell}f(r)+h^{m-\ell}\bigr),
    \qquad 1\leq\ell\leq k,
\]
for $r<R$ sufficiently close to $R$ and $h>0$ sufficiently small.
Here we used $0\leq f(t)\leq f(r)$ on $[r,r+h]$ and boundedness
of $f^{(m)}$ on a fixed neighborhood of $R$.
Taking $h=f(r)^{1/m}$ gives
$|f^{(\ell)}(r)|\leq C_m f(r)^{1-\ell/m}$.
Fa\`a di Bruno's formula, applied on $(0,R)$, therefore implies
\[
    |(f^a)^{(k)}(r)|
    \leq C_{a,k,m}f(r)^{a-k/m}\longrightarrow0
    \qquad\text{as }r\to R.
\]
All derivatives consequently extend by zero across $R$, proving
the claim. Applying it to $f_0$ and $f_1$ shows that
$f_\theta\in C^\infty((0,\infty))$.

Near the origin, where both profiles are positive,
\[
    -\frac{f_\theta'(r)}{f_\theta(r)}
    =(1-\theta)\left(-\frac{f_0'(r)}{f_0(r)}\right)
      +\theta\left(-\frac{f_1'(r)}{f_1(r)}\right)
    \geq\frac{(1-\theta)c_0+\theta c_1}{r}.
\]
The derivative bounds in \textup{(H2)}, together with
Fa\`a di Bruno's formula and Leibniz's rule, also give
\[
    |f_\theta^{(k)}(r)|\leq C_k r^{-k}f_\theta(r),
    \qquad k\geq1,
\]
hence \textup{(H2)} holds.
Finally,
\[
    r^{n+\sigma_\theta-1}\overline\rho_\theta(r)
    =\bigl(r^{n+\sigma_0-1}\overline\rho_0(r)\bigr)^{1-\theta}
     \bigl(r^{n+\sigma_1-1}\overline\rho_1(r)\bigr)^\theta
\]
is almost nonincreasing near the origin. The analogous identity
with $\gamma_0,\gamma_1,\gamma_\theta$ gives an almost
nondecreasing function. This proves \textup{(H3)--(H4)}. \qed

Now, we are in position to establish the complex interpolation identity:
\begin{Teor}
    Let $\rho_0,\rho_1$ be two radial kernels with compact support satisfying hypotheses $(H0)-(H4)$, and $1<p_0,p_1<\infty$. For $0<\theta<1$, define $\rho_\theta(x):=\rho_0(x)^{1-\theta}\rho_1(x)^\theta$ and $1/p_\theta:=(1-\theta)/p_0+\theta/p_1$. Then, with equivalence of norms $$H^{\rho_\theta,p_\theta}(\R^n)=[H^{\rho_0,p_0}(\R^n),H^{\rho_1,p_1}(\R^n)]_\theta,\,0<\theta<1.$$
\end{Teor}
\noindent{}\textbf{Proof:} By Theorem \ref{TEOREMA}, $H^{\rho_j,p_j}(\R^n)=H^{\alpha_{\rho_j}(\cdot),p_j}(\R^n)=F^{\alpha_{\rho_j}(\cdot)}_{p_j,2}$, for $j=0,1$. Now, we define $$\alpha_\theta(t):=\alpha_{\rho_0}(t)^{1-\theta}\alpha_{\rho_1}(t)^\theta,\,t\geq 1.$$ Clearly, $$\alpha_\theta(\langle\xi\rangle)=\alpha_{\rho_0}(\langle\xi\rangle)^{1-\theta}\alpha_{\rho_1}(\langle\xi\rangle)^{\theta}=w_{\rho_0}(\xi)^{1-\theta}w_{\rho_1}(\xi)^\theta.$$
On the other hand, $\alpha_{\rho_\theta}(\langle \xi\rangle)=w_{\rho_\theta}(\xi)$, which by Lemma \ref{lemma-clave}, satisfies \begin{align*}
    w_{\rho_\theta}(\xi)\asymp |\xi|^{1-n}\overline{\rho_\theta}(1/|\xi|)=|\xi|^{(1-\theta)(1-n)}\overline{\rho_0}^{1-\theta}(1/|\xi|)|\xi|^{(\theta)(1-n)}\overline{\rho_1}^{\theta}(1/|\xi|)\asymp w_{\rho_0}(\xi)^{1-\theta}w_{\rho_1}(\xi)^\theta.
\end{align*}
In consequence $\alpha_{\rho_\theta}(\langle\xi\rangle)=w_{\rho_\theta}(\xi)\asymp w_{\rho_0}(\xi)^{1-\theta}w_{\rho_1}(\xi)^\theta=\alpha_{\theta}(\langle\xi\rangle)$. This implies that the spaces $H^{\alpha_\theta(\cdot),p_\theta}(\R^n)$ and $H^{\alpha_{\rho_\theta}(\cdot),p_\theta}(\R^n)=H^{\rho_\theta,p_\theta}(\R^n)$ are equal with equivalence of the norms $\norm{\mathfrak{F}^{-1}(w_{\rho_\theta}\widehat{u})}_{p_\theta}$ and $\norm{\mathfrak{F}^{-1}( w_{\rho_0}^{1-\theta}w_{\rho_1}^\theta\widehat{u})}_{p_\theta}$, which follows from the Mihlin-H\"ormander theorem applied to the multipliers $m_{\rho_0,\rho_1,\theta}$ and $m_{\rho_0,\rho_1,\theta}^{-1}$ where $$m_{\rho_0,\rho_1,\theta}(\xi):=\frac{w_{\rho_\theta}(\xi)}{w_{\rho_0}^{1-\theta}(\xi)w_{\rho_1}^{\theta}(\xi)}.$$

By Theorem \ref{complex}, \begin{align*}H^{\rho_{\theta},p_\theta}(\R^n)&=H^{\alpha_{\rho_\theta}(\cdot),p_\theta}(\R^n)=H^{\alpha_\theta(\cdot),p_\theta}(\R^n)=[H^{\alpha_{\rho_0}(\cdot),p_0}(\R^n),H^{\alpha_{\rho_1}(\cdot),p_1}(\R^n)]_\theta\\
&=[H^{\rho_0,p_0}(\R^n),H^{\rho_1,p_1}(\R^n)]_\theta,\end{align*} as we wanted to prove \qed 

This result in particular implies the embedding $H^{\rho_0,p}(\R^n)\cap H^{\rho_1,p}(\R^n)\xhookrightarrow{}H^{\rho_\theta,p}(\R^n)$, which yields that for every $u\in H^{\rho_0,p}(\R^n)\cap H^{\rho_1,p}(\R^n)$, we have that $$\norm{u}_{H^{\rho_\theta,p}(\R^n)}\leq \norm{u}_{H^{\rho_0,p}(\R^n)}^{1-\theta}\norm{u}_{H^{\rho_1,p}(\R^n)}^\theta.$$

As a particular example we have the following interpolation result:
Let $0<s_0,s_1<1$, $1<p_0,p_1<\infty$, $\lambda>1$ and
$0<\theta<1$. Set
\[
    s_\theta=(1-\theta)s_0+\theta s_1,
    \qquad
    \frac1{p_\theta}=\frac{1-\theta}{p_0}+\frac{\theta}{p_1},
\]
and define
\[
    M_{\lambda,\theta}(\xi)
    :=(\lambda+4\pi^2|\xi|^2)^{s_\theta/2}
      \bigl[\log(\lambda+4\pi^2|\xi|^2)\bigr]^\theta.
\]
Then, with equivalence of norms,
\[
    \begin{aligned}
    &\bigl[H^{s_0,p_0}(\mathbb R^n),
           H_{\log}^{s_1,p_1,\lambda}(\mathbb R^n)\bigr]_\theta\qquad=\left\{u\in\mathcal S'(\mathbb R^n):
       \mathcal F^{-1}(M_{\lambda,\theta}\widehat u)
       \in L^{p_\theta}(\mathbb R^n)\right\},
    \end{aligned}
\]
with the interpolation norm being equivalent to the norm
$\|\mathcal F^{-1}(M_{\lambda,\theta}\widehat u)\|_{L^{p_\theta}}$. Indeed, let $b(t)=\lambda+4\pi^2(t^2-1)$ for $t\ge1$ and take
$\alpha_j(t)=b(t)^{s_j/2}[\log b(t)]^j$, $j=0,1$.
Since $\lambda>1$, both weights belong to $\mathrm{MRO}$.
The first endpoint is $H^{s_0,p_0}$ with equivalent norms, by
Mihlin's theorem applied to the quotient of the two Bessel weights
and its reciprocal. The second endpoint is the logarithmic space
of Section~3.2 by definition. Theorem \ref{TEOREMA} applies and
$\alpha_0(\langle\xi\rangle)^{1-\theta}
 \alpha_1(\langle\xi\rangle)^\theta=M_{\lambda,\theta}(\xi)$.
\subsection{Interpolation of restriction spaces on bounded Lipschitz domains}
In this section we want to properly define a nonlocal space on bounded domains by an equivalent result as the one in Theorem \ref{main} replacing the couple $\left(L^p(\R^n),W^{1,p}(\R^n)\right)$ by $\left(L^p(\Omega),W^{1,p}(\Omega)\right)$. 

To achieve this we would need a common extension operator independent of the parameters of the spaces. There are many extension operators in the literature concerning Sobolev spaces (see the discussion in \cite[Section 2]{ShiYao}); however, the most suitable extension operator is Rychkov's universal linear extension operator, which we shall denote as $\mathcal{E_R}$ introduced by the former in \cite{Rychkov}.
We first observe that for functions $\alpha\in \textbf{RO}$, and $p,q\in [1,\infty)$, $\mathcal{E_R}$ is bounded from $F^{\alpha(\cdot)}_{p,q}(\Omega)$ to $F^{\alpha(\cdot)}_{p,q}(\R^n)$ whenever $\Omega$ has Lipschitz boundary (see \cite[Corollary 3]{AnopMurach}). This allows us to introduce the following spaces:
\begin{Defi}
    Let $\Omega\subset \R^n$, let $\rho$ be a radial kernel satisfying hypotheses $(H0)-(H4)$, and let $1<p<\infty$. Then, we define the nonlocal Sobolev spaces on domains as $$H^{\rho,p}(\Omega):=H^{\rho,p}(\R^n)/\sim,$$ where $f\sim g$ if and only if their restrictions to $\Omega$ coincide. The norm of these spaces is $$\norm{u}_{H^{\rho,p}(\Omega)}=\operatorname{inf}\{\norm{U}_{H^{\rho,p}}: U\in H^{\rho,p}(\R^n), U|_\Omega=u\}.$$
\end{Defi} Observe that if we have a universal extension operator like $\mathcal{E_R}$, and the restriction operator $R_\Omega$, the standard retraction-coretraction allows us to extend Theorem \ref{TEOREMA} for spaces on bounded domains (see \cite[Theorem 18]{AnopMurach}), and hence for the choices $\alpha_\rho$ and $\varphi_{\rho}$, we can define the space $$F^{\rho}_{p,2}(\Omega):=F^{\alpha_\rho(\cdot)}_{p,2}(\Omega)=\langle L^p(\Omega),W^{1,p}(\Omega)\rangle_{\varphi_\rho}.$$ Since by interpolation $\mathcal{E_R}$ is an extension operator for $F^{\rho}_{p,2}(\Omega)$, and $F^{\rho}_{p,2}(\R^n)=H^{\rho,p}(\R^n)$, we conclude that $$H^{\rho,p}(\Omega)=\langle L^p(\Omega),W^{1,p}(\Omega)\rangle_{\varphi_\rho}.$$

\section{Nonlocal Gagliardo spaces for general kernels}
Aside from the Bessel potential spaces $H^{s,p}(\R^n)$, which have attracted a lot of attention from the community of nonlocal PDEs since the works of Shieh and Spector \cite{ShiehSpector2015,ShiehSpector2018}, the major fractional spaces used as an alternative to classical Sobolev spaces are the so-called \textit{Sobolev-Slobodeckij} or \textit{Gagliardo spaces} $W^{s,p}$, defined for an open set $\Omega\subseteq \R^n$ by $$W^{s,p}(\Omega):=\Bigg\{u\in L^p(\Omega): \frac{|u(x)-u(y)|}{|x-y|^{\frac{n}{p}+s}}\in L^p(\Omega\times \Omega)\Bigg\},\,\text{with}\,\norm{u}_{W^{s,p}}:=\left(\norm{u}_p^p+[u]_{W^{s,p}}^p\right)^\frac{1}{p},$$ where $$[u]_{W^{s,p}}:=\left(\int_{\Omega}\int_{\Omega}\frac{|u(x)-u(y)|^p}{|x-y|^{n+sp}}\,dx\,dy\right)^{1/p},$$ is the Gagliardo seminorm. Those spaces are the real interpolation alternative to the complex interpolation ones $H^{s,p}$. We refer to \cite[Theorem~12.5]{Leoni2023} and \cite[Example~1.8]{Lunardi2018} for detailed proofs of this fact, and to \cite{DiNezzaPalatucciValdinoci2012} for a complete discussion of Gagliardo spaces without interpolation theory. In particular, we are interested in their relationship with the Besov scale since $W^{s,p}(\R^n)=B^{s}_{p,p}(\R^n)$ whenever $s\in (0,1)$ and $1\leq p<\infty$. Considering the generalized smoothness scale, we can define the \textit{nonlocal Gagliardo spaces} $W^{\rho,p}(\R^n)$ for a general kernel $\rho$ as
$$W^{\rho,p}(\R^n):=B^{\alpha_\rho(\cdot)}_{p,p}(\R^n),\,1\leq p< \infty.$$ Note that when $p=2$, $$W^{\rho,2}(\R^n)=B^{\alpha_\rho(\cdot)}_{2,2}(\R^n)=F^{\alpha_\rho(\cdot)}_{2,2}(\R^n)=H^{\alpha_\rho(\cdot),2}(\R^n)=H^{\rho,2}(\R^n),$$ which is consistent with the spaces introduced in \cite{BellidoGarciaNiño2026}. These kinds of spaces have been previously considered in \cite{dorer} (see also \cite{carioni}), where they were introduced as follows:
\begin{Defi}[Nonlocal Gagliardo spaces]
Let $1\leq p<\infty$ and let $\rho$ be a radial kernel satisfying hypothesis $(H0)$. For $\Omega\subset \R^n$, we define the \textit{nonlocal Gagliardo spaces} for the kernel $\rho$ $W^{\rho,p}(\Omega)$ as $$W^{\rho,p}(\Omega):=\Bigg\{u\in L^p(\Omega): \frac{u(x)-u(y)}{|x-y|^{1+n(1/p-1)}}\rho(x-y)\in L^p(\Omega\times \Omega)\Bigg\},$$ endowed with the norm $$\norm{u}_{W^{\rho,p}(\Omega)}:=\left(\norm{u}_{L^p(\Omega)}^p+[u]_{\rho,p}^p\right)^{1/p},$$ where the term $$[u]_{\rho,p}:=\left(\int_{\Omega}\int_\Omega \frac{|u(x)-u(y)|^p}{|x-y|^{n-p(n-1)}}\rho(x-y)^p\,dx\,dy\right)^{1/p},$$ is the \textit{nonlocal Gagliardo seminorm}.
\end{Defi}
Standard arguments show that $W^{\rho,p}(\Omega)$ is a Banach space for $1\leq p<\infty$, separable for $1\leq p<\infty$, reflexive for $1<p<\infty$ and Hilbert for $p=2$ (the Hilbertian case will be studied apart in a following section). Those spaces must be understood, as in the Bessel case $H^{\rho,p}(\R^n)$, as a generalization of the so-called fractional Sobolev spaces or Gagliardo spaces $W^{s,p}(\R^n)$ (see \cite{Leoni2023, DiNezzaPalatucciValdinoci2012} and the references therein). In fact, for the choice $\rho=\rho^s=|x|^{-(n+s-1)}$, which produces the Bessel potential spaces in the case $H^{\rho^s,p}$ (see \cite{ShiehSpector2015}), we have that $$[u]_{\rho^s,p}^p=\int_{\Omega}\int_{\Omega}\frac{|u(x)-u(y)|^p}{|x-y|^{n-p(n-1)}}\frac{1}{|x-y|^{np+sp-p}}\,dx\,dy=\int_{\Omega}\int_{\Omega}\frac{|u(x)-u(y)|^p}{|x-y|^{n+sp}}\,dx\,dy,$$ which is the $p$-th power of the seminorm associated to the spaces $W^{s,p}(\Omega)$, and hence $W^{\rho^s,p}(\Omega)=W^{s,p}(\Omega)$ for $1\leq p<\infty$.

In the fractional case, it is known that the spaces $W^{s,p}(\Omega)$ are intermediate with respect to $L^p(\Omega)$ and $W^{1,p}(\Omega)$, and in fact they can be obtained by the real interpolation method under the assumption of Lipschitz boundary for the domain (see \cite[Theorem 12.5]{Leoni2023}, \cite[Example 1.8]{Lunardi2018}).
Our main goal is to prove that under this seminorm, the spaces coincide with the ones from the Besov scale approach. 
\begin{Teor}\label{main2}
    Let $\rho$ satisfy $(H0)-(H4)$ with compact support and $1\leq p<\infty$. We define $$\alpha_\rho(t):=w_\rho(\sqrt{t^2-1}e_1),\,t>1,\,\alpha_\rho(1)=1,$$ and consider $$\varphi_\rho(t):=\begin{cases}
        \alpha_\rho(t^{-1})^{-1},\,0<t<1,\\
        \alpha_\rho(t),\,t\geq 1.
    \end{cases}$$ Then, $\alpha_\rho\in \textbf{MRO}$, $$0<\sigma\leq\sigma_0(\alpha_\rho)\leq\sigma_1(\alpha_\rho)\leq \gamma<1,$$ and the function $\varphi_\rho$ belongs to $\mathfrak{B}$ and satisfies \ref{eqn:condition}. Therefore $$B^{\alpha_\rho(\cdot)}_{p,q}(\R^n)=\left( L^p(\R^n),W^{1,p}(\R^n)\right)_{\varphi_\rho,q},\,1\leq q\leq \infty,$$ and for $1<p<\infty$, $B^{\alpha_\rho(\cdot)}_{p,p}(\R^n)=W^{\rho,p}(\R^n)$ with $$W^{\rho,p}(\R^n)=\{u\in L^p(\R^n): [u]_{\rho,p}<\infty\}.$$
\end{Teor}
\noindent{}\textbf{Proof:} The first part follows from Theorem \ref{main} and Theorem \ref{TEOREMA2}. We only have to prove that the nonlocal Gagliardo seminorm could be written in terms of the generalized real method norm. 
We first recall that by Lemma \ref{lemma-clave} and the discussion on the proof of Theorem \ref{main}, there exists $0<r_0<1$ such that $$|h|^{n-1}\rho(h)\asymp\alpha_\rho(1/|h|),\,0<|h|<r_0,$$ and hence since \begin{align*}[u]^p_{\rho,p}&=\int_{\R^n}\int_{\R^n}\frac{|u(x)-u(y)|^p}{|x-y|^n}\left(|x-y|^{n-1}\rho(x-y)\right)^p\,dy\,dx\\
&=\int_{\R^n}\norm{u(\cdot+h)-u}_p^p\frac{(|h|^{n-1}\rho(h))^p}{|h|^n}\,dh,\end{align*} then \begin{align*}
    [u]^p_{\rho,p}&=\int_{|h|<r_0}\norm{u(\cdot+h)-u}_p^p\frac{(|h|^{n-1}\rho(h))^p}{|h|^n}\,dh+\int_{|h|\geq r_0}\norm{u(\cdot+h)-u}_p^p\frac{(|h|^{n-1}\rho(h))^p}{|h|^n}\,dh\\
    &\leq C\int_{|h|\leq 1}\norm{u(\cdot+h)-u}_p^p\frac{\alpha_\rho(1/|h|)^p}{|h|^n}\,dh+\int_{|h|\geq r_0}\norm{u(\cdot+h)-u}_p^p\frac{(|h|^{n-1}\rho(h))^p}{|h|^n}\,dh.
\end{align*}To bound the second integral, we observe that by hypothesis $(H0)$, we have that for $f_\rho(r)=r^{n-2}\overline{\rho}(r)$, $$\int_{r_0/2}^\infty f_\rho(r)\,dr<\infty,$$ and by $(H1)$, $$\int_{r_0/2}^\infty f_\rho(t)\,dt\geq \int_{r/2}^r f_\rho(t)\,dt\geq \frac{r}{2}f_\rho(r),\,r\geq r_0,$$ which shows that $rf_\rho(r)\leq C$ for every $r\geq r_0$. Then, \begin{align*}
    \int_{r_0}^\infty \frac{(r^{n-1}\overline{\rho}(r))^p}{r}\,dr=\int_{r_0}^\infty (rf_\rho(r))^{p-1}f_\rho(r)\,dr\leq C^{p-1}\int_{r_0}^\infty f_\rho(r)\,dr<\infty.
\end{align*} Finally, changing to polar coordinates, \begin{align*}
    \int_{|h|\geq r_0}\norm{u(\cdot+h)-u}_p^p\frac{(|h|^{n-1}\rho(h))^p}{|h|^n}\,dh\leq 2^p\norm{u}_p^p\int_{|h|\geq r_0}\frac{(|h|^{n-1}\rho(h))^p}{|h|^n}\,dh\leq C\norm{u}_p^p.
\end{align*}

We recall that by Theorem \ref{TEOREMA2}, $$\norm{u}_{B^{\alpha_\rho(\cdot)}_{p,p}}\asymp\left(\int_0^\infty\left(\frac{K(t,u)}{\varphi_\rho(t)}\right)^p\frac{dt}{t}\right)^{1/p},$$ where $$K(t,u)=\inf_{v\in W^{1,p}}\{\norm{u-v}_p+t\norm{v}_{W^{1,p}}\}.$$ Note that for $t\geq 1$, $ct^\sigma\leq \alpha_\rho(t)\leq Ct^\gamma,$ and hence by the definition of $\varphi_\rho$, $$C^{-p}\norm{u}_p^p\int_1^\infty \frac{dt}{t^{1+p\gamma}}\leq\norm{u}_p^p\int_1^\infty \frac{dt}{t\alpha_\rho(t)^p} \leq c^{-p}\norm{u}_p^p\int_1^\infty \frac{dt}{t^{1+p\sigma}},$$ which shows that $$\norm{u}_{B^{\alpha_\rho(\cdot)}_{p,p}}\asymp \norm{u}_p+\left(\int_0^1\left(K(t,u)\alpha_\rho(1/t)\right)^p\frac{dt}{t}\right)^{1/p}.$$ Now, for the first inclusion, we see that for every $v\in W^{1,p}(\R^n)$, $$\norm{u(\cdot+h)-u}_p\leq 2\norm{u-v}_p+|h|\norm{v}_{W^{1,p}},$$ and taking the infimum over $v$, $$\norm{u(\cdot+h)-u}_p\leq 2K(|h|,u),$$ thus \begin{align*}
    \int_{|h|<1}\frac{\norm{u(\cdot+h)-u}_p^p}{|h|^n}\alpha_\rho(1/|h|)^p\,dh&\leq 2^p\int_{|h|<1}\frac{K(|h|,u)^p\alpha_\rho(1/|h|)^p}{|h|^n}\,dh\\
    &= C\int_0^1\frac{K(t,u)^p\alpha_\rho(1/t)^p}{t}\,dt,
\end{align*} where we performed a change to polar coordinates in the last equality. Then, \begin{align*}
    [u]_{\rho,p}^p\lesssim \norm{u}_p^p+\int_0^1\frac{K(t,u)^p\alpha_\rho(1/t)^p}{t}\,dt,
\end{align*} and hence $$\norm{u}_{W^{\rho,p}}\lesssim\norm{u}_p+\left(\int_0^1\frac{K(t,u)^p\alpha_\rho(1/t)^p}{t}\,dt\right)^{1/p}\lesssim \norm{u}_{B^{\alpha_\rho(\cdot)}_{p,p}},$$ which proves the first inclusion.

For the second inclusion, we take $u\in L^p(\R^n)$ such that $[u]_{\rho,p}<\infty$. The mollification estimates in the proof of
\cite[Proposition~4.5, p.~1310]{BrascoSalort2019},
together with Jensen's inequality and the $L^p$-contractivity
of mollification, give, for every $t>0$, $$K(t,u)^p\lesssim t^{-n}\int_{|h|<t}\norm{u(\cdot+h)-u}_p^p\,dh+t^p\norm{u}_p^p.$$
Now we integrate this bound with the weight  \(\alpha_\rho(1/t)^p\,dt/t\), for \(0<t<r_0\), and bound both terms.
For the first term, Fubini's theorem gives us that
\[
\begin{aligned}
&\int_0^{r_0}
\frac{\alpha_\rho(1/t)^p}{t^{n+1}}
\int_{|h|<t}\|u(\cdot-h)-u\|_p^p\,dh\,dt\\
&\quad=
\int_{|h|<r_0}\|u(\cdot-h)-u\|_p^p
\left(\int_{|h|}^{r_0}
\frac{\alpha_\rho(1/t)^p}{t^{n+1}}\,dt\right)dh.
\end{aligned}
\]

For \(|h|\le t<r_0<1\), we already established that
\[
\frac{\alpha_\rho(1/|h|)}{\alpha_\rho(1/t)}
\ge c\left(\frac{t}{|h|}\right)^\sigma
\ge c,
\]
so
\[
\begin{aligned}
\int_{|h|}^{r_0}\frac{\alpha_\rho(1/t)^p}{t^{n+1}}\,dt
&\le C\alpha_\rho(1/|h|)^p
\int_{|h|}^{r_0}\frac{dt}{t^{n+1}}\\
&\le C\frac{\alpha_\rho(1/|h|)^p}{|h|^n}.
\end{aligned}
\]
Now since \(\alpha_\rho(1/|h|)\le C|h|^{n-1}\rho(h)\) for \(|h|<r_0\), we get that
\[
\begin{aligned}
&\int_0^{r_0}
\frac{\alpha_\rho(1/t)^p}{t^{n+1}}
\int_{|h|<t}\|u(\cdot-h)-u\|_p^p\,dh\,dt\\
&\qquad\le
C\int_{|h|<r_0}
\|u(\cdot-h)-u\|_p^p
\frac{\bigl(|h|^{n-1}\rho(h)\bigr)^p}{|h|^n}\,dh
\le C[u]_{\rho,p}^p.
\end{aligned}
\]

For the other term we use the upper bound
\(\alpha_\rho(1/t)\le Ct^{-\gamma}\), where \(\gamma<1\), and we get
\[
\int_0^{r_0}t^p\alpha_\rho(1/t)^p\,\frac{dt}{t}
\le C\int_0^{r_0}t^{p(1-\gamma)-1}\,dt
=\frac{C r_0^{p(1-\gamma)}}{p(1-\gamma)}
<\infty.
\]

Putting everything together yields
\[
\int_0^{r_0}\alpha_\rho(1/t)^pK(t,u)^p\,\frac{dt}{t}
\lesssim[u]_{\rho,p}^p+\|u\|_p^p.
\] The part $$\int_{r_0}^1\alpha_\rho(1/t)^pK(t,u)^p\frac{dt}{t},$$ can be easily bounded by choosing $v=0$ in the definition of the $K$-functional and from the fact that $\alpha_\rho$ is continuous and bounded on $[1,1/r_0],$ which proves the other continuous inclusion. \qed

A direct consequence of this identification is the following embedding result, which generalizes the fractional one (see \cite[Theorem 3.31] {BellidoGarciaI2025}):
\begin{Teor}\label{embeddings-gagliardo-bessel}    Let $\rho$ be a radial kernel with compact support satisfying $(H0)-(H4)$. Then, the following continuous embeddings hold:
    \begin{itemize}
        \item $W^{\rho,p}(\R^n)\xhookrightarrow{}H^{\rho,p}(\R^n)$, $1< p\leq 2$;
        \item $H^{\rho,p}(\R^n)\xhookrightarrow{}W^{\rho,p}(\R^n)$, $2\leq p<\infty$;
        \end{itemize}
       with strict inclusions unless $p=2$.
\end{Teor}
\noindent{}\textbf{Proof:} The result follows directly from Theorem \ref{EmbeddingsNL}, Theorem \ref{main} and Theorem \ref{main2}. \qed

We also have the inclusion of the nonlocal Gagliardo spaces between fractional Gagliardo spaces as it was partially proved in \cite[Proposition 4.5]{dorer}.
\begin{Teor}
    Let $\rho$ be a general kernel with compact support satisfying $(H0)-(H4)$ and $1<p<\infty$. Then, $$W^{\gamma,p}(\R^n)\xhookrightarrow{}W^{\rho,p}(\R^n)\xhookrightarrow{}W^{\sigma,p}(\R^n).$$
\end{Teor}
\noindent{}\textbf{Proof:} By Theorem \ref{TEOREMA2} and \ref{main2}, we know that $$W^{\rho,p}(\R^n)=\left(W^{\sigma,p}(\R^n),W^{\gamma,p}(\R^n)\right)_{\varphi_{\rho,\sigma,\gamma},p},$$ where $$\varphi_{\rho,\sigma,\gamma}(t):=t^{-\frac{\sigma}{\gamma-\sigma}}\varphi_\rho\left(t^{1/(\gamma-\sigma)}\right),$$ and 
$\sigma$ and $\gamma$ are the fractional parameters in hypotheses $(H3)$ and $(H4)$, respectively. Then, since $W^{\gamma,p}(\R^n)\xhookrightarrow{}W^{\sigma,p}(\R^n)$ by \cite[Proposition 2.1]{DiNezzaPalatucciValdinoci2012}, we have that $$W^{\gamma,p}(\R^n)\cap W^{\sigma,p}(\R^n)=W^{\gamma,p}(\R^n)\xhookrightarrow{}W^{\rho,p}(\R^n)\xhookrightarrow{}W^{\gamma,p}(\R^n)+W^{\sigma,p}(\R^n)=W^{\sigma,p}(\R^n),$$ as we wanted to prove. \qed

As it happened in the Bessel case, Theorem \ref{TEOREMA2} also holds for bounded domains $\Omega$ with Lipschitz boundary (see \cite[Theorem 7]{AnopMurach}), and hence we can define the nonlocal Gagliardo spaces on domains as $W^{\rho,p}(\Omega)=B^{\alpha_\rho(\cdot)}_{p,p}(\Omega)$, and by Theorem \ref{main2}, $$W^{\rho,p}(\Omega)=(L^p(\Omega),W^{1,p}(\Omega))_{\varphi_\rho,p},\,1\leq p<\infty,$$ whenever $\Omega$ has Lipschitz boundary.
\subsection{Some properties of nonlocal Gagliardo spaces}
In this subsection we will study several properties of the spaces $W^{\rho,p}$ without interpolation techniques, adapting the classical proofs from the fractional setting and paying special attention to the hypothesis over the kernels.
The first result is the embedding between spaces when the kernels are related by their behavior near zero:
\begin{Prop}\label{GagliardoAnid}
    Let $\Omega\subset \R^n$ be open, let $1\leq p<\infty$ and let $\rho_0,\rho_1$ be two radial kernels satisfying $(H0)$ and such that $$\liminf_{t\to 0}\frac{\overline{\rho_1}(t)}{\overline{\rho_0}(t)}>0.$$ Then, if some of the following holds:
    \begin{itemize}
        \item $p=1$;
        \item $p>1$ and $\rho_0$ satisfies $(H1)$;
    \end{itemize} then the continuous embedding $$W^{\rho_1,p}(\Omega)\xhookrightarrow{}W^{\rho_0,p}(\Omega),$$holds.
\end{Prop}\textbf{Proof:} By the hypothesis $\liminf_{t\to 0}\frac{\overline{\rho_1}(t)}{\overline{\rho_0}(t)}>0,$ there exist $c>0$ and $r_0>0$ such that for each $0<|h|<r_0$, $\rho_1(h)\geq c\rho_0(h)$. Hence, \begin{align}\label{eqn:1}
    &\int_\Omega\int_{\Omega\cap\{|x-y|<r_0\}}\frac{|u(x)-u(y)|^p}{|x-y|^{n-p(n-1)}}\rho_0(x-y)^p\,dx\,dy \\
    &\leq c^{-p}\int_\Omega\int_{\Omega\cap\{|x-y|<r_0\}}\frac{|u(x)-u(y)|^p}{|x-y|^{n-p(n-1)}}\rho_1(x-y)^p\,dx\,dy\\
    &\leq c^{-p}[u]_{\rho_1,p}^p.
\end{align} On the other hand, since $|u(x)-u(y)|^p\leq 2^{p-1}\left(|u(x)|^p+|u(y)|^p\right)$,\begin{align}\label{eqn:2}
    &\int_\Omega\int_{\Omega\cap\{|x-y|\geq r_0\}}\frac{|u(x)-u(y)|^p}{|x-y|^{n-p(n-1)}}\rho_0(x-y)^p\,dx\,dy\\
    &\leq 2^{p-1}\int_\Omega\int_{\Omega\cap\{|x-y|\geq r_0\}}\frac{|u(x)|^p+|u(y)|^p}{|x-y|^{n-p(n-1)}}\rho_0(x-y)^p\,dx\,dy\\
    &\leq 2^p\int_\Omega\left(\int_{|h|\geq r_0}\frac{\rho_0(h)^p}{|h|^{n-p(n-1)}}dh\right)|u(x)|^p\,dx.
\end{align} Let us write $$I_{r_0}:=\int_{|h|\geq r_0}\frac{\rho_0(h)^p}{|h|^{n-p(n-1)}}dh.$$ First, if we have that $p=1$, the integral reduces to $$\int_{|h|\geq r_0}\frac{\rho_0(h)}{|h|}dh,$$ which is finite by hypothesis $(H0)$ under the kernel $\rho_0$. If $p>1$, we require $\rho_0$ under $(H1)$ in order to have that $f_{\rho_0}(r)=r^{n-2}\overline{\rho_0}(r)$ is decreasing. By $(H0)$, the integral $$\int_{r_0/2}^\infty f_{\rho_0}(t)\,dt:=M_{r_0},$$ is finite. Then, by $(H1)$, for every $r\geq r_0$, $$\frac{r}{2}f_{\rho_0}(r)\leq \int_{r/2}^rf_{\rho_0}(t)\,dt\leq M_{r_0},$$ and hence $r^{n-1}\overline{\rho_0}(r)=rf_{\rho_0}(r)\leq 2 M_{r_0}<\infty.$ Then, for any $|h|\geq r_0$, we have that $$|h|^{p(n-1)-n}\rho_0(h)^p=\left(|h|^{n-1}\rho_0(h)\right)^{p-1}\frac{\rho_0(h)}{|h|}\leq (2M_{r_0})^{p-1}\frac{\rho_0(h)}{|h|},$$ and hence $$I_{r_0}\leq (2M_{r_0})^{p-1}\int_{B_{r_0}(0)^c}\frac{\rho_0(h)}{|h|}\,dh,$$ which is finite by $(H0)$ again. In both cases, $I_{r_0}<\infty$. Using this in \ref{eqn:2}, we get that $$\int_{\Omega}\int_{\Omega\cap \{|x-y|\geq r_0\}}\frac{|u(x)-u(y)|^p}{|x-y|^{n-p(n-1)}}\rho_0(x-y)^p\,dx\,dy\leq 2^pI_{r_0}\norm{u}_{L^p(\Omega)}^p.$$ Combining \ref{eqn:1} and \ref{eqn:2}, \begin{align*}
    \norm{u}_{W^{\rho_0,p}}^p=\norm{u}_{L^p(\Omega)}^p+[u]_{\rho_0,p}^p\leq \norm{u}_{L^p(\Omega)}^p+c^{-p}[u]_{\rho_1,p}^p+2^pI_{r_0}\norm{u}_{L^p(\Omega)}^p\leq C\norm{u}_{W^{\rho_1,p}}^p,
\end{align*} for $$C:=\max\{1+2^pI_{r_0},c^{-p}\},$$ and hence the continuous embedding $$W^{\rho_1,p}(\Omega)\xhookrightarrow{}W^{\rho_0,p}(\Omega),$$ holds. \qed

\begin{Obse}
In the case of fractional kernels $\rho_0=\rho^{s_0}, \rho_1=\rho^{s_1}$, with $0<s_0<s_1<1$, we have that $\liminf_{t\to 0}\overline{\rho_1}(t)/\overline{\rho_0}(t)=\infty$, and hence the previous result gives us that $$W^{\rho^{s_1},p}(\Omega)\xhookrightarrow{}W^{\rho^{s_0},p}(\Omega),\,1\leq p<\infty,$$ which recovers the well-known fact of the nesting of the scale $W^{s,p}(\Omega)$ on the fractional parameter (see \cite[Proposition 2.1]{DiNezzaPalatucciValdinoci2012}).
\end{Obse}
  
The next result shows, like in the fractional case, that the local space $W^{1,p}(\Omega)$ embeds continuously into $W^{\rho,p}(\Omega)$ under the assumption of $\Omega$ being an extension domain and the existence of a universal extension operator, for example for special Lipschitz domains (see \cite{Rychkov}).
\begin{Prop}\label{SoboGagliardo}
    Let $\Omega$ be a special Lipschitz domain, $1\leq p<\infty$ and let $\rho$ be a radial kernel satisfying $(H0)-(H1)$. Then, $$W^{1,p}(\Omega)\xhookrightarrow{}W^{\rho,p}(\Omega).$$
\end{Prop}
\noindent{}\textbf{Proof:} Let $u\in W^{1,p}(\Omega)$. By hypothesis, there exists a function $U:\R^n\to \R$ such that $U\in W^{1,p}(\R^n)$, $U|_\Omega=u$ and $\norm{U}_{W^{1,p}(\R^n)}\leq C\norm{u}_{W^{1,p}(\Omega)}$ for a constant only depending on $\Omega$. We split the seminorm again and bound each one separately.
First, \begin{align*}
    \int_\Omega\int_{\Omega\cap \{|x-y|<1\}}\frac{|u(x)-u(y)|^p}{|x-y|^{n-p(n-1)}}\rho(x-y)^p\,dx\,dy\leq\int_\Omega\int_{B_1(0)}\frac{|U(x+h)-U(x)|^p}{|h|^{n-p(n-1)}}\rho(h)^p\,dh\,dx.
\end{align*}
By the fundamental theorem of calculus and H\"older's inequality, we know that $$|U(x+h)-U(x)|^p\leq |h|^p\int_0^1|\nabla U(x+th)|^p\,dt,$$ and hence \begin{align*}
    &\int_\Omega\int_{\Omega\cap \{|x-y|<1\}}\frac{|u(x)-u(y)|^p}{|x-y|^{n-p(n-1)}}\rho(x-y)^p\,dx\,dy \\
    &\leq \int_\Omega\int_{B_1(0)}\left(\int_0^1|\nabla U(x+th)|^p\,dt\right)\frac{\rho(h)^p}{|h|^{n-np}}\,dh\,dx\\
    &\leq \int_{\R^n}\int_{B_1(0)}\left(\int_0^1|\nabla U(x+th)|^p\,dt\right)\frac{\rho(h)^p}{|h|^{n-np}}\,dh\,dx\\
    & \leq \norm{\nabla U}_{L^p(\R^n)}^p\int_{B_1(0)}\left(\int_0^1\,dt\right)\frac{\rho(h)^p}{|h|^{n-np}}\,dh=\norm{\nabla U}_{L^p(\R^n)}^p\int_{B_1(0)}\frac{\rho(h)^p}{|h|^{n-np}}\,dh\\
    &\leq C\norm{u}_{W^{1,p}(\Omega)}^p\int_{B_1(0)}\frac{\rho(h)^p}{|h|^{n-np}}\,dh=C\omega_{n-1} \norm{u}_{W^{1,p}(\Omega)}^p\int_0^1 r^{np-1}\overline{\rho}(r)^p\,dr.
\end{align*}
Since $\rho$ verifies $(H0)$, the integral $$M:=\int_0^1r^{n-1}\overline{\rho}(r)\,dr,$$ is finite. Moreover, since it verifies $(H1)$, for every $0<t<r<1$, $$r^{n-2}\overline{\rho}(r)\leq t^{n-2}\overline{\rho}(t),$$ so \begin{align*}
    2M=2\int_0^1t^{n-1}\overline{\rho}(t)\,dt\geq 2\int_0^r t^{n-1}\overline{\rho}(t)\,dt\geq 2\int_0^r  tr^{n-2}\overline{\rho}(r)\,dt=r^n\overline{\rho}(r),
\end{align*}
so $$\int_0^1r^{np-1}\overline{\rho}(r)^p\,dr=\int_0^1(r^n\overline{\rho}(r))^{p-1}r^{n-1}\overline{\rho}(r)\,dr\leq (2M)^{p-1}\int_0^1 r^{n-1}\overline{\rho}(r)\,dr=2^{p-1}M^p.$$ Putting this together, we get that $$\int_\Omega\int_{\Omega\cap \{|x-y|<1\}}\frac{|u(x)-u(y)|^p}{|x-y|^{n-p(n-1)}}\rho(x-y)^p\,dx\,dy\leq C\omega_{n-1}2^{p-1}M^p\norm{u}_{W^{1,p}(\Omega)}^p.$$ For the other part of the seminorm, arguing as in the previous result, we got that $$\int_{\Omega}\int_{\Omega\cap \{|x-y|\geq 1\}}\frac{|u(x)-u(y)|^p}{|x-y|^{n-p(n-1)}}\rho(x-y)^p\,dx\,dy\leq 2^pI_{1}\norm{u}_{L^p(\Omega)}^p,$$ which is finite under the hypothesis $(H1)$ over the kernel. Hence, there exists a positive constant $\widetilde{C}>0$ such that $$[u]_{\rho,p}^p\leq \widetilde{C}\norm{u}_{W^{1,p}(\Omega)}^p,$$ thus the desired embedding holds. \qed

\begin{Obse}
    If the domain is not Lipschitz, we can construct a function in $W^{1,p}(\Omega)$ that does not belong to $W^{\rho,p}(\Omega)$. In particular, for the choice $\rho=\rho^s$, such an example is given in \cite[Example 9.1]{DiNezzaPalatucciValdinoci2012}.
\end{Obse}

By a standard method of approximation by convolution (see \cite[Theorem 6.62]{Leoni2023}), which translates to our setting without further assumptions, we get the following density result:
\begin{Prop}
    Let $\rho$ satisfy $(H0)$ and $1\leq p<\infty$. Then, it holds that $C^\infty\cap W^{\rho,p}(\R^n)$ is dense in $W^{\rho,p}(\R^n)$.
\end{Prop}

By the well-known
mollification and cut-off arguments (see \cite[Theorem 6.66]{Leoni2023}), we have the density of $C_c^\infty$ functions:
\begin{Prop}
     Let $\rho$ satisfy $(H0)$ and $(H1)$, and $1\leq p<\infty$. Then, it holds that $C_c^\infty(\R^n)$ is dense in $W^{\rho,p}(\R^n)$.
\end{Prop}

A relevant remark which connects our spaces with other families of nonlocal Gagliardo spaces is the condition $$M_{\rho,p}:=\int_{\R^n}\min\{1,|h|^p\}\frac{\rho(h)^p}{|h|^{n-p(n-1)}}\,dh<\infty,\,1\leq p<\infty.$$ Note that when $p=1$ this reduces to the hypothesis $(H0)$. In the proof of Proposition \ref{GagliardoAnid}, we used that under the hypothesis $(H0)-(H1)$, the kernel $\rho_0$ satisfies $$\int_{|h|\geq r_0}\frac
{\rho_0(h)^p}{|h|^{n-p(n-1)}}\,dh<\infty,$$ and in Proposition \ref{SoboGagliardo}, that the same hypothesis implies $$\int_{|h|<1}|h|^p\frac{\rho(h)^p}{|h|^{n-p(n-1)}}\,dh<\infty.$$ Both estimates combined are equivalent to the condition $M_{\rho,p}<\infty$, which is known as the \textit{$p$-L\'evy integrability condition} (see \cite{Foghem2025}), and hence we have that under the hypothesis $(H0)-(H1)$, the kernel $|h|^{-(n-p(n-1))}\rho(h)^p$ is a symmetric $p$-L\'evy measure.

In \cite{Foghem2025}, Gagliardo type spaces modeled on seminorms given by symmetric $p$-L\'evy measures are studied, and we can directly translate several structural results to our setting under the hypothesis $(H0)-(H1)$, of particular relevance being the fact that the condition $M_{\rho,p}<\infty$ is equivalent to the embedding $W^{1,p}(\R^n)\xhookrightarrow{}W^{\rho,p}(\R^n)$ being continuous \cite[Theorem 4.1]{Foghem2025}, and the fact that $W^{\rho,p}(\R^n)\not =\{0\}$ if and only if $\rho$ satisfies $M_{\rho,p}<\infty$ \cite[Theorem 4.3]{Foghem2025}, which as we recalled is ensured by $(H0)-(H1)$.
\subsection{The Hilbertian case}
In this subsection we focus on the case $p=2$, which as in the fractional and classical cases is of particular interest due to the broader class of techniques available such as Fourier characterizations. 

As we stated previously, we know that under the hypothesis $(H0)-(H4)$, the spaces $W^{\rho,2}(\R^n)$ and $H^{\rho,2}(\R^n)$ coincide with equivalence of norms because of the equivalence $B^{\alpha_\rho(\cdot)}_{2,2}(\R^n)=F^{\alpha_\rho(\cdot)}_{2,2}(\R^n)$ and the identifications with those spaces. Our objective now is to prove this fact without the use of interpolation theory, analogously as \cite[Proposition 4.1]{DiNezzaPalatucciValdinoci2012}. We recall the potential characterization of $H^{\rho,2}(\R^n)$ as $$H^{\rho,2}(\R^n)=\{u\in \mathcal{S}'(\R^n): \left(1+|\lambda_\rho (\xi)|^2\right)^{1/2}\widehat{u}(\xi)\in L^2(\R^n)\},$$ which is straightforward by Plancherel's theorem and the fact that $\widehat{D_\rho u}(\xi)=\lambda_\rho (\xi)\widehat{u}(\xi).$

We start by characterizing the spectral definition of the nonlocal Gagliardo space:
\begin{lema}
    Let $\rho$ satisfy $(H0)$. Then, for every $u\in W^{\rho,2}(\R^n)$ we have that $$\norm{u}_2^2+[u]_{\rho,2}^2=\int_{\R^n}\left(1+m_\rho(\xi)\right)|\widehat{u}(\xi)|^2\,d\xi,$$ where $$m_\rho(\xi)=2\int_{\R^n}\frac{\rho(h)^2}{|h|^{2-n}}\left(1-\cos(2\pi \xi\cdot h)\right)\,dh,\,\xi\in \R^n.$$
\end{lema}
\noindent{}\textbf{Proof:} Let $x\in \R^n$ be fixed and $u\in W^{\rho,2}(\R^n)$. By the change of variables $x-y=h$ and Plancherel's formula: \begin{align*}
[u]_{\rho,2}^2&=\int_{\R^n}\int_{\R^n}\frac{|u(x)-u(y)|^2\rho(x-y)^2}{|x-y|^{2-n}}\,dy\,dx=\int_{\R^n}\int_{\R^n}\frac{|u(x+h)-u(x)|^2\rho(h)^2}{|h|^{2-n}}\,dh\,dx\\
    &=\int_{\R^n}\frac{\rho(h)^2}{|h|^{2-n}}\left(\int_{\R^n}|u(x+h)-u(x)|^2\,dx\right)\,dh\\
    &=2\int_{\R^n}\frac{\rho(h)^2}{|h|^{2-n}}\left(\int_{\R^n}|\widehat{u}(\xi)|^2\left(1-\cos(2\pi \xi\cdot h)\right)\,d\xi\right)\,dh\\
&=2\int_{\R^n}\left(\int_{\R^n}\frac{\rho(h)^2}{|h|^{2-n}}\left(1-\cos(2\pi \xi\cdot h)\right)\,dh\right)|\widehat{u}(\xi)|^2\,d\xi,
\end{align*}
as we wanted to prove. \qed

Now the equivalence of the spaces relies on the equivalence of the symbols $m_\rho$ and $|\lambda_\rho|^2$:
\begin{lema}
    Let $\rho$ satisfy $(H0)-(H4)$ with compact support. Then, $$1+m_{\rho}(\xi)\asymp 1+|\lambda_\rho(\xi)|^2.$$
\end{lema}
\noindent{}\textbf{Proof:} We first estimate $m_\rho$ from above. Recall that
\[
m_\rho(\xi)
=
2\int_{\mathbb{R}^n}
\bigl(1-\cos(2\pi\xi\cdot z)\bigr)
|z|^{n-2}\rho(z)^2\,dz.
\]
Fix $\xi\neq 0$ and split the integral into the regions
$|z|<|\xi|^{-1}$ and $|z|\geq|\xi|^{-1}$.

For the first integral, the inequality
\[
1-\cos t=2\sin^2(t/2)\leq \frac{t^2}{2}
\]
gives
\[
2\bigl(1-\cos(2\pi\xi\cdot z)\bigr)
\leq 4\pi^2|\xi\cdot z|^2
\leq 4\pi^2|\xi|^2|z|^2.
\]
Therefore,
\begin{align*}
&2\int_{|z|<|\xi|^{-1}}
\bigl(1-\cos(2\pi\xi\cdot z)\bigr)
|z|^{n-2}\rho(z)^2\,dz\qquad\leq
4\pi^2|\xi|^2
\int_{|z|<|\xi|^{-1}}|z|^n\rho(z)^2\,dz.
\end{align*}
For the second integral, using $1-\cos t\leq 2$, we obtain
\begin{align*}
&2\int_{|z|\geq|\xi|^{-1}}
\bigl(1-\cos(2\pi\xi\cdot z)\bigr)
|z|^{n-2}\rho(z)^2\,dz\qquad\leq
4\int_{|z|\geq|\xi|^{-1}}|z|^{n-2}\rho(z)^2\,dz.
\end{align*}
Changing to polar coordinates yields
\[
\int_{|z|<|\xi|^{-1}}|z|^n\rho(z)^2\,dz
=
\omega_{n-1}
\int_0^{|\xi|^{-1}}
r^{2n-1}\overline{\rho}(r)^2\,dr,
\]
and
\[
\int_{|z|\geq|\xi|^{-1}}|z|^{n-2}\rho(z)^2\,dz
=
\omega_{n-1}
\int_{|\xi|^{-1}}^\infty
r^{2n-3}\overline{\rho}(r)^2\,dr.
\]
Recalling that
\[
f_\rho(r)=r^{n-2}\overline{\rho}(r),
\]
we have
\[
r^{2n-1}\overline{\rho}(r)^2=r^3f_\rho(r)^2,
\qquad
r^{2n-3}\overline{\rho}(r)^2=rf_\rho(r)^2,
\]
and combining the previous estimates, we conclude that
\[
m_\rho(\xi)
\leq C\left(
|\xi|^2\int_0^{|\xi|^{-1}}r^3f_\rho(r)^2\,dr
+
\int_{|\xi|^{-1}}^\infty rf_\rho(r)^2\,dr
\right),
\]
where $C>0$ depends only on the dimension.

We estimate the first integral.
By (H4), the function
\[
r\longmapsto r^{n+\gamma-1}\overline{\rho}(r)
=r^{1+\gamma}f_\rho(r)
\]
is almost nondecreasing on $(0,\varepsilon)$, so there
exists $C>0$ such that
\[
r^{1+\gamma}f_\rho(r)
\leq C t^{1+\gamma}f_\rho(t),
\qquad 0<r\leq t<\varepsilon.
\]
For $|\xi|>\varepsilon^{-1}$, we may take $t=|\xi|^{-1}$
and obtain
\[
f_\rho(r)
\leq C |\xi|^{-1-\gamma}r^{-1-\gamma}
f_\rho(|\xi|^{-1}),
\qquad 0<r\leq|\xi|^{-1}.
\]
Consequently,
\begin{align*}
|\xi|^2\int_0^{|\xi|^{-1}}r^3f_\rho(r)^2\,dr
&\leq
C^2|\xi|^2|\xi|^{-2-2\gamma}
f_\rho(|\xi|^{-1})^2
\int_0^{|\xi|^{-1}}r^{3-2-2\gamma}\,dr\\
&=
C^2|\xi|^{-2\gamma}
f_\rho(|\xi|^{-1})^2
\int_0^{|\xi|^{-1}}r^{1-2\gamma}\,dr.
\end{align*}
Since $\gamma<1$, the last integral is finite and equals
\[
\int_0^{|\xi|^{-1}}r^{1-2\gamma}\,dr
=
\frac{(|\xi|^{-1})^{2-2\gamma}}{2-2\gamma}
=
\frac{|\xi|^{-2+2\gamma}}{2-2\gamma}.
\]
Therefore,
\[
|\xi|^2\int_0^{|\xi|^{-1}}r^3f_\rho(r)^2\,dr
\leq
\frac{C^2}{2-2\gamma}
|\xi|^{-2}f_\rho(|\xi|^{-1})^2,
\]
and recalling the definition of $f_\rho$,
\[
|\xi|^2\int_0^{|\xi|^{-1}}r^3f_\rho(r)^2\,dr
\leq
C_\gamma|\xi|^{2-2n}\overline{\rho}(1/|\xi|)^2.
\]
We next estimate the second integral. Fix $r_0\in(0,\varepsilon)$
and choose $R>r_0$ such that $\operatorname{supp}\rho\subset B_R(0)$.
For $|\xi|>r_0^{-1}$, hypothesis (H3) gives
\[
r^{1+\sigma}f_\rho(r)
\leq C|\xi|^{-1-\sigma}f_\rho(|\xi|^{-1}),
\qquad |\xi|^{-1}\leq r\leq r_0.
\]
For $r_0\leq r\leq R$, the monotonicity of $f_\rho$ in $(H1)$
and the previous inequality evaluated at $r=r_0$ yield
\begin{align*}
r^{1+\sigma}f_\rho(r)
&\leq R^{1+\sigma}f_\rho(r_0)=\left(\frac{R}{r_0}\right)^{1+\sigma}
r_0^{1+\sigma}f_\rho(r_0)\leq C\left(\frac{R}{r_0}\right)^{1+\sigma}
|\xi|^{-1-\sigma}f_\rho(|\xi|^{-1}).
\end{align*}
Since $f_\rho(r)=0$ for $r\geq R$, we conclude that
\[
f_\rho(r)
\leq C|\xi|^{-1-\sigma}r^{-1-\sigma}
f_\rho(|\xi|^{-1}),
\qquad r\geq|\xi|^{-1},
\]
where $C>0$ is independent of $\xi$ and $r$. Consequently, since $\sigma>0$,
\begin{align*}
\int_{|\xi|^{-1}}^\infty rf_\rho(r)^2\,dr
&\leq C^2|\xi|^{-2-2\sigma}f_\rho(|\xi|^{-1})^2
\int_{|\xi|^{-1}}^\infty r^{-1-2\sigma}\,dr=\frac{C^2}{2\sigma}
|\xi|^{-2}f_\rho(|\xi|^{-1})^2,
\end{align*} so
we obtain
\[
\int_{|\xi|^{-1}}^\infty rf_\rho(r)^2\,dr
\leq C_\sigma|\xi|^{2-2n}\overline{\rho}(|\xi|^{-1})^2.
\]
Combining both estimates we get that 
\[
m_\rho(\xi)
\leq C|\xi|^{2-2n}\overline{\rho}(|\xi|^{-1})^2,
\qquad |\xi|>r_0^{-1},
\]
where $C>0$ depends on $n$ and $\rho$, but not on $\xi$.

Now, by \cite[Lemma 4.3]{BMS} $$\widehat{Q}_\rho(\xi)\geq C|\xi|^{-n}\overline{\rho}(1/|\xi|),\,|\xi|\geq \frac{2}{\varepsilon},$$ which implies that $$|\lambda_\rho(\xi)|^2\geq 4\pi^2 C|\xi|^{2-2n}\overline{\rho}(1/|\xi|)^2,\,|\xi|\geq \frac{2}{\varepsilon},$$ and hence for $|\xi|\geq \max\{2/\varepsilon,r_0^{-1}\}$, we have that $$m_\rho(\xi)\leq C |\lambda_\rho(\xi)|^2.$$ 

We next prove the lower bound. Fix $\xi\in\mathbb R^n\setminus\{0\}$. By concavity of the sine function on $[0,\pi/2]$,
\[
\sin s\ge \frac{2s}{\pi},
\qquad 0\le s\le \frac{\pi}{2},
\]
and consequently,
\[
1-\cos t
=
2\sin^2\left(\frac{|t|}{2}\right)
\ge \frac{2}{\pi^2}t^2,
\qquad |t|\le\pi.
\]
If $|z|<\frac12|\xi|^{-1}$, then $|2\pi\xi\cdot z|<\pi$, and hence
\[
2\bigl(1-\cos(2\pi\xi\cdot z)\bigr)
\ge 16(\xi\cdot z)^2.
\]
Since the integrand defining $m_\rho$ is nonnegative, we obtain
\[
\begin{aligned}
m_\rho(\xi)
&=
2\int_{\mathbb R^n}
\bigl(1-\cos(2\pi\xi\cdot z)\bigr)
|z|^{n-2}\rho(z)^2\,dz\\
&\ge 2\int_{|z|<\frac12|\xi|^{-1}}
\bigl(1-\cos(2\pi\xi\cdot z)\bigr)
|z|^{n-2}\rho(z)^2\,dz\\
&\ge
16\int_{|z|<\frac12|\xi|^{-1}}
(\xi\cdot z)^2|z|^{n-2}\rho(z)^2\,dz.
\end{aligned}
\]

Writing $z=r\theta$ in polar coordinates and using the radiality of $\rho$, we find
\[
m_\rho(\xi)
\ge
16
\left(
\int_{\mathbb S^{n-1}}(\xi\cdot\theta)^2\,d\theta
\right)
\left(
\int_0^{\frac12|\xi|^{-1}}
r^{2n-1}\bar\rho(r)^2\,dr
\right).
\]
By rotational invariance and symmetry,
\[
\int_{\mathbb S^{n-1}}(\xi\cdot\theta)^2\,d\theta
=
|\xi|^2\int_{\mathbb S^{n-1}}\theta_1^2\,d\theta
=
\frac{\omega_{n-1}}{n}|\xi|^2.
\]
Indeed, the integrals of $\theta_j^2$ are equal for $j=1,\ldots,n$, and
$\sum_{j=1}^n\theta_j^2=1$ on $\mathbb S^{n-1}$.
Since $f_\rho(r)=r^{n-2}\bar\rho(r)$, it follows that
\[
m_\rho(\xi)
\ge
\frac{16\omega_{n-1}}{n}
|\xi|^2
\int_0^{\frac12|\xi|^{-1}}r^3f_\rho(r)^2\,dr.
\]
By (H1), $f_\rho$ is nonincreasing, hence
\[
f_\rho(r)\ge f_\rho(|\xi|^{-1}),
\qquad 0<r<\frac12|\xi|^{-1},
\]
and therefore
\[
\begin{aligned}
\int_0^{\frac12|\xi|^{-1}}r^3f_\rho(r)^2\,dr
&\ge
f_\rho(|\xi|^{-1})^2
\int_0^{\frac12|\xi|^{-1}}r^3\,dr=
\frac{|\xi|^{-4}}{64}f_\rho(|\xi|^{-1})^2.
\end{aligned}
\]
Combining these estimates yields
\[
m_\rho(\xi)
\ge
\frac{\omega_{n-1}}{4n}
|\xi|^{-2}f_\rho(|\xi|^{-1})^2.
\]
Finally, using
\[
f_\rho(|\xi|^{-1})
=
|\xi|^{2-n}\bar\rho(|\xi|^{-1}),
\]
we conclude that
\[
m_\rho(\xi)
\ge
\frac{\omega_{n-1}}{4n}
|\xi|^{2-2n}\bar\rho(|\xi|^{-1})^2,
\qquad \xi\ne0.
\]
Finally, by \cite[Lemma 7.1]{BMS}, there exists a positive constant $C$ such that $$\widehat{Q}_\rho(\xi)\leq C|\xi|^{-n}\overline{\rho}(1/|\xi|),\,|\xi|\geq \varepsilon^{-1},$$ and hence $$|\lambda_\rho(\xi)|^2\leq C |\xi|^{2-2n}\overline{\rho}(1/|\xi|)^2,\,|\xi|\geq \varepsilon^{-1},$$ so for $|\xi|\geq \varepsilon^{-1}$, $$m_\rho(\xi)\geq C|\lambda_\rho(\xi)|^2.$$ Combining both the upper and lower estimates on the symbol $m_\rho$, we conclude that for $|\xi|\geq \max\{2/\varepsilon,r_0^{-1}\}:=R_0$, $m_\rho(\xi) \asymp|\lambda_\rho(\xi)|^2$.
It remains to consider the region $|\xi|\le R_0$.
To compare the weights $1+m_\rho(\xi)$ and $1+|\lambda_\rho(\xi)|^2$ in this region, it suffices to show that both symbols are bounded there.
We first consider $m_\rho$. Using $1-\cos t\le t^2/2$, we obtain
\[
\begin{aligned}
m_\rho(\xi)
&=
2\int_{\mathbb R^n}
\bigl(1-\cos(2\pi\xi\cdot z)\bigr)
|z|^{n-2}\rho(z)^2\,dz\le
4\pi^2\int_{\mathbb R^n}
(\xi\cdot z)^2|z|^{n-2}\rho(z)^2\,dz\\
&\le
4\pi^2|\xi|^2
\int_{\mathbb R^n}|z|^n\rho(z)^2\,dz.
\end{aligned}
\]
The last integral is finite. Indeed, passing to polar coordinates gives
\[
\int_{\mathbb R^n}|z|^n\rho(z)^2\,dz
=
\omega_{n-1}
\int_0^\infty r^{2n-1}\bar\rho(r)^2\,dr
=
\omega_{n-1}
\int_0^\infty r^3f_\rho(r)^2\,dr.
\]
Since $0<r_0<\varepsilon$, assumption (H4) yields
\[
f_\rho(r)\le Cr^{-1-\gamma},
\qquad 0<r<r_0,
\]
where $C>0$ is independent of $r$. Hence
\[
\int_0^{r_0}r^3f_\rho(r)^2\,dr
\le
C^2\int_0^{r_0}r^{1-2\gamma}\,dr
=
\frac{C^2r_0^{2-2\gamma}}{2-2\gamma}
<\infty.
\]
Moreover, by (H1), $f_\rho(r)\le f_\rho(r_0)$ for $r\ge r_0$. Since $\rho$ has compact support, it follows that
\[
\int_{r_0}^\infty r^3f_\rho(r)^2\,dr<\infty,
\]
and consequently,
\[
\sup_{|\xi|\le R_0}m_\rho(\xi)
\le
4\pi^2R_0^2
\int_{\mathbb R^n}|z|^n\rho(z)^2\,dz
<\infty.
\]
Now, since $\rho$ has compact support, $Q_\rho \in L^1(\R^n)$, hence $$|\lambda_\rho(\xi)|^2=4\pi^2|\xi|^2\widehat{Q}_\rho(\xi)^2\leq 4\pi^2R_0^2 \norm{Q_\rho}_1^2<\infty,\,0<|\xi|\leq R_0,$$ which concludes the proof. \qed 

As a direct consequence we obtain the identification of the spaces $H^{\rho,2}$ and $W^{\rho,2}$ with equivalence of norms:
\begin{Coro}
    Let $\rho$ satisfy $(H0)-(H4)$ with compact support. Then the spaces $H^{\rho,2}(\R^n)$ and $W^{\rho,2}(\R^n)$ are equal with equivalence of norms.
\end{Coro}
\subsection{A broader class of kernels}
In the recent work \cite{Brusca}, the authors study spaces of non-local Gagliardo-type associated
with integral kernels under suitable conditions. In particular, they introduce the spaces $$W^{K,p}(\Omega):=\{u\in L^p(\Omega):[u]_{K,p}<\infty\},\,1\leq p<\infty,$$ where $$[u]_{K,p}:=\left(\int_\Omega\int_\Omega|u(x)-u(y)|^pK(x-y)\,dx\,dy\right)^{1/p}.$$ Here, $K$ is a measurable function from $\R^n\to [0,+\infty]$ that is not zero up to negligible sets satisfying the following conditions:
\begin{itemize}
    \item[(Nint)] $K\not \in L^1(\R^n)$;
    \item[(Far)] $K\in L^1(\R^n\setminus B_R)$, for every $R>0$;
    \item[(Dec)] There exists $c_0\in(0,1]$ such that for every $|x|\leq|y|$, $K(x)\geq c_0K(y)$;
    \item[(Dou)] There exist positive constants $D$ and $C_D$ such that for every $|x|\leq D$, $K(x)\leq C_DK(2x)$;
    \item[(Nts$_p$)] $\int_{\R^n}\min\{1,|z|^p\}K(z)\,dz<\infty$.
\end{itemize}
We refer to \cite{Brusca} and the references therein for a discussion on these hypotheses. Under these hypotheses, existence of extension operators, compactness results and Poincar\'e inequalities are established. Our main goal in this subsection will be to relate our spaces for general kernels $\rho$ with the ones considered in \cite{Brusca}, and hence translate directly their results to our setting under the suitable hypothesis over $\rho.$

Let $1\leq p<\infty$. We denote $K_{\rho,p}(z):=|z|^{p(n-1)-n}\rho(z)^p$.
We begin with the first property of nonintegrability. It is hypothesis $(H3)$, which implies that our kernels lie above a purely fractional kernel around zero.
\begin{lema}
    Let $\rho$ be a radial kernel satisfying $(H0)$ and $(H3).$ Then, $K_{\rho,p}$ satisfies (Nint) for every $1\leq p<\infty$.
\end{lema}
\noindent{}\textbf{Proof:} By $(H3)$, there exists a positive constant $C$ such that $\overline{\rho}(t)\geq Ct^{-(n+\sigma-1)},$ for every $t\in (0,\varepsilon)$, where $\varepsilon$ is the number in hypothesis $(H0)$. Indeed, fix $t_0\in (0,\varepsilon)$. By $(H3)$, for every $0<t<t_0$, $$t^{n+\sigma-1}\overline{\rho}(t)\geq c(t_0)t_0^{n+\sigma-1}\overline{\rho}(t_0):=C(t_0),$$ and hence $\overline{\rho}(t)\geq C(t_0)t^{-(n+\sigma-1)}$. Then we can just take $C$ as $C:=\operatorname{sup}_{t_0\in (0,\varepsilon)}C(t_0)$, which is finite by $(H0)$. Now, we have that for every $0<|z|<\varepsilon$, \begin{align*}
    K_{\rho,p}(z)=|z|^{p(n-1)-n}\rho(z)^p\geq C|z|^{p(n-1)-n}|z|^{-p(n+\sigma-1)}=C|z|^{-(n+p\sigma)},
\end{align*}
 thus \begin{align*}
     \int_{\R^n}K_{\rho,p}(z)\,dz\geq \int_{B_\varepsilon(0)}K_{\rho,p}(z)\,dz\geq C\int_{B_\varepsilon(0)}|z|^{-n-p\sigma}\,dz=C\omega_{n-1}\int_0^\varepsilon r^{-1-p\sigma}\,dr=\infty,
 \end{align*} and the result follows. \qed

Now, we look at the (Far) hypothesis.
\begin{lema}
     Let $\rho$ be a radial kernel satisfying $(H0)$ and $(H1).$ Then, $K_{\rho,p}$ satisfies (Far) for every $1\leq p<\infty$.
\end{lema}
\noindent{}\textbf{Proof:} Let $R>0$. Without loss of generality we can assume that $R/2\geq 1$. If not, suppose that $R/2<1$. Hypothesis $(H0)$ implies that $$\int_0^1tf_\rho(t)\,dt<\infty,\,\int_1^\infty f_\rho(t)\,dt<\infty.$$ Then, $$\int_{R/2}^\infty f_\rho(t)\,dt=\int_{R/2}^1f_\rho(t)\,dt+\int_1^\infty f_\rho(t)\,dt,$$ and $$\int_{R/2}^1f_\rho(t)\,dt\leq \frac{2}{R}\int_{R/2}^1tf_\rho(t)\,dt\leq \frac{2}{R}\int_0^1tf_\rho(t)\,dt<\infty,$$ by the monotonicity in $(H1)$. Hence, we can restrict ourselves to the case $R/2\geq 1$. Now, for $t\geq R$, since $f_\rho(s)\geq f_\rho(t)$ for $s\in [t/2,t]$, we have that $$\frac{t}{2}f_\rho(t)\leq\int_{t/2}^tf_\rho(s)\,ds\leq \int_{R/2}^\infty f_\rho(s)\,ds:=A_R<\infty,$$ thus $tf_\rho(t)\leq 2A_R$ for every $t\in (0,\infty)$. Finally, \begin{align*}
    \int_{|z|\geq R}K_{\rho,p}(z)\,dz&=\omega_{n-1}\int_{R}^\infty t^{p(n-1)-1}\overline{\rho}(t)^p\,dt=\omega_{n-1}\int_R^\infty \left(tf_\rho(t)\right)^{p-1}f_\rho(t)\,dt \\
    &\leq (2A_R)^{p-1}\omega_{n-1}\int_R^\infty f_\rho(t)\,dt\leq 2^{p-1}A_R^p\omega_{n-1}<\infty,
\end{align*}
for every $R\geq 1.$ \qed

With these hypotheses, the following compactness result follows directly from \cite[Theorem 5.1]{Brusca}:
\begin{Teor}[Compactness on the whole space]
Let $\rho$ be a kernel satisfying $(H0), (H1)$ and $(H3)$. Let $p\in [1,\infty)$ and let $(u_k)_k\subset W^{\rho,p}(\R^n)$ be a sequence such that $$\operatorname{sup}_k\norm{u_k}_{W^{\rho,p}}<\infty.$$ Then, there exists a subsequence $(u_{k_j})_j\subset(u_k)_k$ and $u\in W^{\rho,p}$ such that $u_{k_j}\to u$ in $L^p_{loc}(\R^n)$ as $j\to\infty.$
    
\end{Teor}
As discussed in the previous subsections, from the proofs of Proposition \ref{GagliardoAnid} and Proposition \ref{SoboGagliardo}, we have that $M_{\rho,p}<\infty$, and hence the following:
\begin{lema}
     Let $\rho$ be a radial kernel satisfying $(H0)$ and $(H1).$ Then, $K_{\rho,p}$ satisfies (Nts$_p$) for every $1\leq p<\infty$.
\end{lema}

We now look at the (Dou) condition:
\begin{lema}
    Let $\rho$ be a radial kernel satisfying $(H0)$ and $(H2).$ Then, $K_{\rho,p}$ satisfies (Dou) for every $1\leq p<\infty$.
\end{lema}
\noindent{}\textbf{Proof:}
By $(H2)$, there exists a positive constant $C_1$ such that $-f_\rho' (t)\leq C_1\frac{f_\rho(t)}{t},\,t\in (0,\varepsilon)$. Let $D:=\varepsilon/4$ and $0<r\leq D$. Then, the interval $[r,2r]\subset(0,\varepsilon)$ and \begin{align*}
    \log\frac{f_\rho(r)}{f_\rho(2r)}=-\int_r^{2r}\frac{f'_\rho(t)}{f_\rho(t)}\,dt\leq C_1\int_r^{2r}\frac{dt}{t}=C_1\log 2,
\end{align*}
which implies that $f_\rho(r)\leq 2^{C_1}f_{\rho}(2r)$ for every $r\in (0,D]$. Now, let us denote by $\overline{K}_{\rho,p}$ the radial representation of $K_{\rho,p}$, i.e.,$r^{p(n-1)-n}\overline{\rho}(r)^p$. We have that \begin{align*}
    \overline{K}_{\rho,p}(r)=r^{p-n}f_{\rho}(r)^p\leq r^{p-n}2^{pC_1}f_\rho(2r)^p=r^{p-n}2^{pC_1}(2r)^{n-p}\overline{K}_{\rho,p}(2r)=2^{n+p(C_1-1)}\overline{K}_{\rho,p}(2r),
\end{align*} and the result follows choosing $C_D:=2^{n+p(C_1-1)}$. \qed

Finally, we investigate the condition (Dec):
\begin{lema}
    Let $\rho$ be a kernel satisfying:
    \begin{itemize}
        \item If $p\leq n$, $(H0)$ and $(H1)$;
        \item If $p>n$, $(H0), (H1), (H3)$ and compact support;
    \end{itemize} Then, $K_{\rho,p}$ satisfies (Dec).
\end{lema}
\noindent{}\textbf{Proof:}
First, we consider the case $p\leq n$. In that case, for $r\leq t$, we have that since $n-p\geq 0$, $t^{p-n}\leq r^{p-n}$, and by $(H1)$, $f_\rho(t)\leq f_\rho (r)$. Hence, $$\overline{K}_{\rho,p}(t)=t^{p-n}f_\rho(t)^{p}\leq r^{p-n}f_\rho(r)^p=\overline{K}_{\rho,p}(r),$$ and hence (Dec) follows for $c_0=1$.

For $p>n$, we will consider that $\operatorname{supp}\rho \subset B_R(0)$ for some $R>0$. Let $0<a<\min\{\varepsilon,R\}$. First, let $0<r\leq t\leq a$. By $(H3)$, $$t^{n+\sigma-1}\overline{\rho}(t)\leq Cr^{n+\sigma-1}\overline{\rho}(r),$$ for some positive constant $C$, that we can choose greater than $1$. Then, \begin{align*}
    \overline{K}_{\rho,p}(t)\leq C^p t^{p(n-1)-n}\left(\frac{r}{t}\right)^{p(n+\sigma-1)}\overline{\rho}(r)^p=C^p\left(\frac{r}{t}\right)^{n+p\sigma}\overline{K}_{\rho,p}(r)\leq C^p\overline{K}_{\rho,p}(r).
\end{align*}
For $a\leq r\le t\le R$, we have that $$\overline{K}_{\rho,p}(t)\leq t^{p-n}f_\rho(r)^p=\left(\frac{t}{r}\right)^{p-n}\overline{K}_{\rho,p}(r)\leq \left(\frac{R}{a}\right)^{p-n}\overline{K}_{\rho,p}(r).$$

If $r<a<t\leq R$, combining the previous estimates, \begin{align*}
    \overline{K}_{\rho,p}(t)\leq \left(\frac{R}{a}\right)^{p-n}\overline{K}_{\rho,p}(a)\leq C^p\left(\frac{R}{a}\right)^{p-n}\overline{K}_{\rho,p}(r).
\end{align*} Finally, if $t>R,$ then $\overline{K}_{\rho,p}(t)=0$, so we have covered all the possible cases and hence the result follows choosing $c_0:=C^{-p}\left(\frac{a}{R}\right)^{p-n}\leq 1$. \qed

In view of the previous result, if we consider a kernel $\rho$ with compact support and satisfying $(H0)-(H3)$, we have that the spaces $W^{\rho,p}(\Omega)$ satisfy the hypothesis of \cite[Theorem 4.5]{Brusca}, and hence the following result follows directly:
\begin{Teor}
    Let $\Omega$ be a bounded, connected open set with Lipschitz boundary, $1\leq p<\infty$ and consider a kernel $\rho$ with compact support and satisfying $(H0)-(H3)$. Then:
    \begin{itemize}
        \item For every $u\in W^{\rho,p}(\Omega)$, there exists $U\in W^{\rho,p}(\R^n)$ such that $U|_\Omega=u$ and $\norm{U}_{W^{\rho,p}(\R^n)}\leq C\norm{u}_{W^{\rho,p}(\Omega)}$, for some positive constant $C$ depending on $\Omega,p$ and $\rho$;
        \item $W^{\rho,p}(\Omega)$ is compactly embedded into $L^p(\Omega)$;
        \item There exists a positive constant $C$ such that $$\left(\int_\Omega|u-u_\Omega|^p\,dx\right)^{1/p}\leq C[u]_{\rho,p},$$ for every $u\in W^{\rho,p}(\Omega),$ where $u_\Omega=|\Omega|^{-1}\int_\Omega u\,dx.$
    \end{itemize}
\end{Teor}
\noindent{}\textbf{Proof:} It follows from the previous lemmas and \cite[Theorem 4.5, Corollary 5.2, Proposition 5.6]{Brusca}. \qed
 \addcontentsline{toc}{section}{References}
	% Incluir la bibliografía
	\bibliographystyle{plain} 
	
	\addcontentsline{toc}{section}{Statements \& Declarations}
	\subsection*{Funding}
    This work is supported by {\it Agencia Estatal de Investigación} (Spain) through grant PID2023-151823NB-I00 and {\it Junta de Comunidades de Castilla-La Mancha} (Spain) through grant SBPLY/23/180225/000023. G.G.-S. is supported by a Doctoral Fellowship by Universidad de Castilla-La Mancha.  \text{2024-UNIVERS-12844-404}.
\section*{Acknowledgements}
The author would like to thank his advisor José Carlos Bellido for the useful suggestions on a preliminary
version of the manuscript and Leon Happ, Carolin Kreisbeck and Carlos Mora-Corral for several fruitful discussions on the topic.

	\subsection*{Conflicts of interest}
	The author declares that there are no conflicts of interest regarding the publication of this paper.

\subsection*{Generative AI declaration}The author proudly declares that no generative
AI tools were used in the development of the work. All the mathematical ideas, the theoretical derivation of results and all the proofs are the result of the human mind. 

\par\bigskip
\begingroup
\footnotesize
\setlength{\parindent}{1em}

\indent
\textsc{Guillermo García-Sáez.
Departamento de Matemáticas,
Universidad de Castilla-La Mancha,
Campus Universitario s/n, 13071, Ciudad Real, Spain.}\par

\nobreak
\indent
\textit{Email address}:
\texttt{\detokenize{guillermo.garciasaez@uclm.es}}\par

\endgroup
\end{document}